\documentclass{amsart} 

\usepackage{amsmath, amsthm}
\usepackage{amssymb, latexsym}
\usepackage{amsfonts} 
\usepackage[dvipsnames]{xcolor}
 
\usepackage{threeparttable}

    \usepackage{lineno}
\newtheorem{Theorem}[equation]{Theorem}
\newtheorem{theorem}[equation]{Theorem}
\newtheorem{Proposition}[equation]{Proposition}

\newtheorem{lemma}[equation]{Lemma}
\newtheorem{Lemma}[equation]{Lemma}

\newtheorem{question}[equation]{Question}
\newtheorem{questions}[equation]{Questions}

\newtheorem{Example}[equation]{Example}

\newtheorem{Remark}[equation]{Remark}
\newtheorem{Remarks}[equation]{Remarks}

\theoremstyle{definition}
\newtheorem{definition}[equation]{Definition}
\newtheorem{notation}[equation]{Notation}

\newcommand{\exref}[1]{Ex\-am\-ple \ref{#1}}
\newcommand{\thmref}[1]{\rm Theorem~\ref{#1}}

\newcommand{\secref}[1]{Sec\-tion \ref{#1}}
\newcommand{\lemref}[1]{Lem\-ma \ref{#1}}
\newcommand{\propref}[1]{\rm Proposition~\ref{#1}}

\newcommand{\remref}[1]{Re\-mark \ref{#1}}

\newcommand{\beql}[1]{\begin{equation}\label{#1}}
\newcommand{\eeq} {\end{equation}}

    \font\Aaa=msam10

\font\Bbb=msbm10

\newcommand\F{\hbox{\Bbb F}}

\numberwithin{equation}{section}

\newcommand\g{\gamma }
\newcommand\G{\Gamma }
\newcommand\D{{ \Delta }}

        \DeclareMathOperator{\rad}{rad}

        \def\O{{\rm O}}

        \def\SL{{\rm SL}}

        \def\Sp{{\rm Sp}}

        \def\qed{\hbox{~~\Aaa\char'003}}

        \def\Sz{{\rm Sz}}
               
        \def\<{{\langle}}
        \def\>{{\rangle}}

 \def\a{\alpha}
\def\b{\beta}

\def\div{ \kern-.5pt\hbox{\big |} }
\def\ndiv{ {\not\kern-.5pt\hbox{\big |}\,} }
\def\ndivv{ {\not\kern+1.5pt\hbox{$\mid$}\,} }

 \def\rad{{\rm rad}}
\def\dim{{\rm dim}}

\def\col{\colon\!}

\def\B{^2\kern-.8pt B}
\def\G{^2\kern-.8pt G}
\def\EH{^2\kern-.8pt\hat  E}
\def\E{^2\kern-.8pt E}
\def\D{^3\kern-1pt D}
\def\FF{^2\kern-.8pt F}

\newdimen\refcodesize
\newbox\seriesbox
\def\para#1{\medskip\noindent{\bf #1}~}

\def\proof{\noindent {\bf Proof.~}}

 \font\Bbb=msbm10

 \font\Bbb=msbm10

\def\F{\hbox{\Bbb F}}
\font\Bbbsmall=msbm7

\def\Zsmall{\hbox{\Bbbsmall Z}}

\def\I{{\mathfrak I}}

\DeclareRobustCommand{\SkipTocEntry}[4]{}
 
\usepackage{lineno}
 
\begin{document} 

\title[maximal symplectic  partial  spreads]
{On maximal symplectic  partial  spreads 
 }

       \author{W. M. Kantor
      }
       \address{U. of Oregon,
       Eugene, OR 97403
        \      and  \
       Northeastern U., Boston, MA 02115}
       \email{kantor@uoregon.edu}

\begin{abstract}
\vspace{-4pt}
New types of maximal symplectic  partial spreads are constructed.

 \end{abstract}

\maketitle
\vspace{-20pt}
   \section{Introduction
   }
      \label{Introduction}

Since  very few papers  concern  maximal symplectic partial spreads in dimension $>4$
\cite{Grassl},    this paper will focus on those   dimensions.  
The largest and most obvious type of maximal partial spread of a 
$2n$-dimensional symplectic $\F_q$-space is a spread, of size $q^n+1$, which we will not consider here. 
(However, there are relatively few known types of symplectic spreads, see \cite{Ka3} for a survey as of 2012.) 

On the other hand, 
when $n$ is even
Grassl   \cite{Grassl}   initially conjectured that the smallest possible size
of a maximal symplectic partial spread
 is $q^{n/2} +1$, and he provided examples
of this size for all even $q$  and $n$.  However, 
when $2n=8$ the conjecture is not correct \cite{Grassl}.  
Families of counterexamples using   Suzuki-Tits ovoids 
are   in \secref{Grassl's Conjecture}.
It still seems plausible that  Grassl's conjecture may  be correct 
if   $2n>8$     or  if $q$ is  odd.  
Thus far all   counterexamples  to this conjecture  have size greater than
$q^{n/2} /2$.

     
Most of our  examples   are based on  standard properties of orthogonal and symplectic spaces, involving either orthogonal spreads or
the standard method for  obtaining them
(Sections~\ref{Using orthogonal  spreads},
 \ref{More maximal  symplectic partial spreads}
 and  \ref{Projections}),
or partial $\O^+(8,q)$-ovoids and triality (Section~\ref{$8$-dimensional  partial spreads}).
   Approximately  half of this paper is concerned with spaces of dimension 4 or 8,  where we can  use points as crutches:
the Klein correspondence
in dimension 4  \cite[p.~196]{Taylor} and triality  in dimension 8 
\cite{Ti}
 turn  sets of points into sets of subspaces
 (of dimension 2 or 4).
In dimension $>4$ our results are 
summarized in Table~\ref{Dimension at least 6}; 
the pairs of dimensions of the form $4n,$ $4n-2$ 
arise  from orthogonal partial spreads and 
are explained in~\secref{Projections}.

 Maximal  symplectic partial spreads have a straightforward
 use in  Quantum Physics  for finding 
 sets of mutually unbiased  bases \cite{MBGW,Grassl} (e.g., by plugging into 
   \cite[Eqs.~(3.2)~or~(3.4)]{Ka3}  in  order~to obtain sets of complex vectors).
  
 There are tables of computer-generated sizes of maximal symplectic partial spreads  
 in  $\F_q^{2n}$~for very small $n$ and $q$   \cite{CDFS,Grassl}.
 A~few of these are special cases of  constructions given here. 
 However,   since
 these tables contain  integer
 intervals  that consist  of  sizes
 of these partial spreads, it is clear that 
 new types  of construction techniques are needed in all dimensions.
 

 \vspace{-3pt}
  \section{Background}
    \label{Background}

  The letter $q$ will always denote a prime power, while $n,m,k,s$ and $i$ will be integers.

\


~\!\!\!\!  \!\!  \!\!  
\begin{threeparttable} [h] 
\renewcommand{\TPTminimum}{\linewidth}
 \caption{\!Maximal symplectic partial spreads:\!   
 dimension $\ge6$
 over~$\F_q$}
\label{Dimension at least 6}
\begin{tabular}{|c|c|c|c|c|c|c|l|l|c|c|c|c|c|c|l|}
\hline 
dimensions			&parity of  $q$&  
Size       &
   Restrictions    &\!Theorems\!
\raisebox{2.4ex}{\hspace{-1pt}}\raisebox{-1.1ex}{\hspace{-1pt}}
 \\
\hline \hline
$ 4m$& arbitrary           & $\!q^{2m}-q^m+(2,q-1) \!$
 &   
 &  \ref{using transversals}       
\raisebox{2.7ex}{\hspace{-1pt}}\raisebox{-1.4ex}{\hspace{-1pt}}
\\ 
\hline 
$\! 4mk , 4mk-2\!$& even  &$ q^{2 m k-k}+1 $
&$ m > (k+1)/2  $  &
\ref{nk}, \ref{project k}
\raisebox{2.7ex}{\hspace{-1pt}}\raisebox{-1.4ex}{\hspace{-1pt}}
    \\
    \hline 
$ 4k$,\,$ 4k-2$& even  &$q^k+1\,$\tnote{a}
&   &
\ref{Grassl example}, \ref{project symplectic1}
\raisebox{2.7ex}{\hspace{-1pt}}\raisebox{-1.4ex}{\hspace{-1pt}}
     \\
\hline 
$ 4k$& even  &$2q^k+1$ 
&   &
\ref{Grassl example}
\raisebox{2.7ex}{\hspace{-1pt}}\raisebox{-1.4ex}{\hspace{-1pt}}
    \\
\hline 
      \mbox{$ 8$ and $6 $}& even   &$q^3 -q^2+1 $
&   $q \ge4$   &
\ref{triality to spread},  \ref{Sp6 corollary}
 \raisebox{2.7ex}{\hspace{-1pt}}\raisebox{-1.4ex}{\hspace{-1pt}}
    \\ 
        \hline 
      \mbox{$ 8$ and $6 $}& even    &$  n_s$\tnote{b}
&$1\le s\le q/ 5 $  &
\ref{orthovoids}, \ref{Sp6 corollary}
\raisebox{2.7ex}{\hspace{-1pt}}\raisebox{-1.4ex}{\hspace{-1pt}}
    \\
    \hline 
      \mbox{$ 8$ and $6 $}& even    & $n_4-1$
& $q\ge16 $  &
\ref{orthovoids}, \ref{Sp6 corollary}
\raisebox{2.7ex}{\hspace{-1pt}}\raisebox{-1.4ex}{\hspace{-1pt}}
    \\ 
        \hline 
      \mbox{$ 8$ and $6 $}& even  &$q^2+1$
&   &
\ref{partial ovoid from spread}, \ref{Sp6 corollary}
 \raisebox{2.7ex}{\hspace{-1pt}}\raisebox{-1.4ex}{\hspace{-1pt}}
    \\
    \hline 
      \mbox{$ 8$ and $6 $}& even  &$2q^2+1$
&   &
\ref{2q2+1}, \ref{Sp6 corollary}
\raisebox{2.7ex}{\hspace{-1pt}}\raisebox{-1.4ex}{\hspace{-1pt}}
    \\
        \hline 
      \mbox{$ 8$ and $6 $}& even  &$q^2+q+1\ $
& $q=2^{2e+1}>2$  &
\ref{q2+q+1}, \ref{Sp6 corollary}
\raisebox{2.7ex}{\hspace{-1pt}}\raisebox{-1.4ex}{\hspace{-1pt}}
    \\
        \hline 
      \mbox{$ 8$ and $6 $}& even  &$q^2-q+1 $
&  $q=2^{2e+1}>2$ &
  \ref{easiest Grassl}, \ref{Sp6 corollary}
\raisebox{2.7ex}{\hspace{-1pt}}\raisebox{-1.4ex}{\hspace{-1pt}}
    \\
    \hline 
   \mbox{$ 8$ and $6 $}& even  & $ q^2-sq+2s-1 $
&  $q=2^{2e+1}>2$ &
\ref{smaller Grassl}, \ref{Sp6 corollary}
 \vspace{-4pt}
\raisebox{2.7ex}{\hspace{-1pt}}\raisebox{-1.4ex}{\hspace{-1pt}}
    \\ 
      &    & 
&  $1 <  s\le  2^{e} -1$ &
\raisebox{2.7ex}{\hspace{-1pt}}\raisebox{-1.4ex}{\hspace{-1pt}}
    \\
    \hline 
$ 6$& arbitrary  &$q^3-q^{2}+1$
&  
&
\ref{group}
\raisebox{2.7ex}{\hspace{-1pt}}\raisebox{-1.4ex}{\hspace{-1pt}}
    \\
    \hline 
\end{tabular}
\begin{tablenotes}
\item [a] {\scriptsize This corresponds   to the excluded possibility $m=1$ in dimensions  $4mk,$   $4mk-2$.}
\item [b] {\scriptsize$n_s=q^3-sq^2+(s-1)(q+2)+\binom{s}{2}(q-2)+1 $}
\vspace{10pt}
\end{tablenotes} 
 \end{threeparttable}

 See  \cite{Taylor} for the standard properties of the 
 symplectic and 
 orthogonal vector spaces  used here.  We name  geometries  using their  isometry groups.  
We will be concerned with singular vectors  and totally singular (t.s.) 
subspaces of  orthogonal spaces, and totally isotropic (t.i.) subspaces of   symplectic spaces.  
A subspace of an orthogonal space is 
\emph{anisotropic} if it contains no nonzero singular vector -- and hence has dimension $\le2$.
In characteristic 2, an orthogonal   vector space   is also 
a  symplectic space, and   t.s. subspaces are also t.i.  subspaces.

 \para{Types of maximal  t.s.~subspaces. \ }
The $n$-dimensional t.s.~subspaces of an $\O^+(2n,q)$-space
 are of two types,~with two such subspaces of  the same type if and only if  their intersection has  dimension
$\equiv n$ (mod 2).   
Each  t.s.~$n -1$-space is contained in one member of each type.  
Since we will be concerned with subspaces intersecting in 0,
 $n$ will be even.

A triality map for an $\O^+(8,q)$-space \cite{Ti}  permutes the t.s. subspaces, sending  singular points to a type of t.s.~4-spaces and   non-perpendicular 
pairs of points to pairs of  4-spaces having zero intersection.
   
   \para{Partial ovoids and partial spreads. }
A \emph{partial ovoid} of an 
 orthogonal  space
 is a set $\Omega$ of   t.s. points such that each
  maximal~t.s.~subspace 
contains  at most one point in the set;  $\Omega$  is an  \emph{ovoid} if it meets every such subspace.   A \emph{partial spread} in a $2n$-dimensional vector space $V$
is a set 
$\Sigma$ of 
$n$-spaces any two of which have only 0 in common;
$\Sigma$  is a \emph{spread} if every   vector is 
in a member of $\Sigma$.~If $V$ is a $2n$-dimensional symplectic or orthogonal  vector space,~a 
 \emph{symplectic or orthogonal   partial spread}
$\Sigma$ is a partial spread consisting~of t.i.~or t.s.~$n$-spaces;
$\Sigma$ is a  \emph{symplectic or  orthogonal    spread} if every vector or 
 every  singular  vector is in a member of $\Sigma$.  This note concerns 
 \emph{maximal}  symplectic  or  orthogonal   partial spreads: maximal with respect to inclusion.  In  some situations we will even obtain  
 symplectic  maximal partial spreads:  maximal partial spreads that happen to be symplectic.

Two symplectic partial spreads are 
\emph{equivalent} if there is a semilinear automorphism of the symplectic geometry sending one partial spread to the other.  
 If $\Sigma $ is a set of subspaces of an $\Sp(2n,q)$-space, then
$\Sp(2n,q)_\Sigma$ is its set-stabilizer in the symplectic group   $\Sp(2n,q)$.    
There are similar definitions for orthogonal spaces and 
for the  automorphism group  of a symplectic or orthogonal partial spread.


\section{Maximal partial $\Sp(4m,q)$-spreads}
\label{Maximal partial $Sp(4m,q)$-spreads}

Our most general result is the following

 \Theorem
 \label{using transversals}
For any $q  $  and $m\ge1,$  an $\Sp(4m,q)$-space has a maximal
 symplectic partial spread of size $q^{2m}-q^m+(2,q-1)$.
 
\rm\medskip

We begin with notation.  
 Let $F=\F_{q^{2 m}}\supset E=\F_{q^{m}} \supset K=\F_q$,
  with trace map  $T\col F\to K$, so that
    $T(xy)$ is a nondegenerate symmetric $K$-bilinear form on $F$.
  By dimensions, $\{x\in F\mid T(x E)=0\} = \theta  E$ 
  for some $\theta\in  F$.  
  
  Equip  the $K$-space $V =F^2 $ with the nondegenerate alternating $K$-bilinear form
  $f\big ((x,y),(x',y') \big):=T(xy')-T(x'y)$.
  
    Let $\Sigma$ be the desarguesian symplectic spread of  
 $V  $  consisting of the t.i. $2$-spaces $[x=0]$
  and $[y=ax]$ for $a\in F$.  
  Let $Z_\star$ be the t.i. $2m$-space $(E,\theta E) = E\oplus \theta E$ (t.i. since 
  $T(E\theta E)=0$).
  
 Let  $\Sigma_\star \subset  \Sigma$  
  consist of the   members of $\Sigma$ met nontrivially by $Z_\star$  (namely, the   $2m$-spaces 
 $[x=0]$ and $[y=a\theta x]$ for $a\in E).$  We need information concerning   some   transversals of $\Sigma_\star$:
 
\lemma
\label{one transversal}
There are exactly $(2,q-1)$  t.i.~$2m$-spaces of the $\Sp(4m,q)$-space $V$  that
meet each member of $\Sigma_\star$ in an $m$-space.  If there are two such subspaces then they intersect in $0$.

\rm\medskip
\proof      
 If $Z$ is such a subspace let $Z\cap [y=0]= (U,0)$
 and $Z\cap [x=0]= (0,W)$ for $m$-dimensional $K$-subspaces $U$
 and $W$ of $F$.
 Since $Z=(U,0) + (0,W)$ is t.i.~we have $T(UW)=0$.
 
 Since $Z\cap [y=a\theta x] $ (for $a\in E$) consists of the vectors $(u,a\theta u)$ with $u\in U$, 
we see that
 $W=\theta U$   (using $a=1$)  and $W$   is closed under multiplication by elements of $E$.  Then
 $W$ is an  
  $E$-subspace of $F$.   Let $U=\a E,$ $ \a\in F^*$, so that 
  $W=\theta \a E$.
 Then $0=T(UW ) = T( \a \theta \a E ) $, 
 so that  $\a^2\theta \in  \theta E$.
 Thus, there are $(2,q-1)$ choices for   the coset $\a F^*\in F^*/E^*$,
 and  hence also  $(2,q-1)$ choices for 
 $Z=(U,  W)=(\a E,  \theta \a E)$.  
 
 This argument reverses:  if the coset $\a E^*$ has order at most 2, then 
 $(\a E,  \theta \a E)$  is a t.i. $2m$-space that 
 meets each member of $\Sigma_\star$ in an $m$-space.
 
 Finally, if there are two subspaces $( E,  \theta  E)$ and $ (\a E,  \theta \a E)$, 
 then $\a\notin E$ and these have intersection $0$.
 \qed
 
  \rm\medskip
{\noindent \bf   Proof of 
Theorem~\ref{using transversals}. }    Let $\Sigma $  and  $\Sigma_\star$  be as above.
By the lemma, there are t.i. $2m$-spaces $Z$ (if $q$ is even) or
$Z,Z'$ (if $q$ is odd)  such that 
$\Sigma_\star$ is the set of elements of $\Sigma$ met nontrivially by either of these
$2m$-spaces. Then
\vspace{-2pt}
$$\Sigma^\bullet:=
  \begin{cases} 
  (\Sigma - \Sigma_\star )\cup \{Z\} &\mbox{if $q$ is even}   
  \vspace{2pt}
    \\
 (\Sigma - \Sigma_\star )\cup \{Z,Z'\} 	&\mbox{if $q$ is odd}  
\vspace{-2pt}
 \end{cases}
 $$
 is a symplectic partial spread  of size
 $ q^{2m}-q^m+(2,q-1)$.
 
\emph{Maximality}:  Suppose that $X$ is a t.i.~$2m$-space meeting each member of $\Sigma^\bullet$ in zero.
Since $\Sigma $ is a spread,  the set  $\Sigma_X$ of members of 
$\Sigma$ meeting $X$ nontrivially must be contained in $\Sigma_\star$.
If \  $(\star) $$\,\Sigma_X = \Sigma_\star $ \emph{and $\dim X\cap Y=0$ or $m$ for each}
 $Y \in \Sigma, $ then $X=Z$ or $Z'$  by \lemref{one transversal},
 which contradicts the fact that  $X\notin \Sigma^\bullet$.

We count
in order to  prove $(\star) $.   Let $a_i$ be the number   of $Y\in \Sigma$
such that  $\dim X\cap Y=i$, where $1\le i\le 2m-1$.  Since
the intersections $X\cap Y$ produce a partition of $X-\{0\}$,  
\vspace{-2pt}
$$
\begin{array}{rll}
\displaystyle \sum_1^{2m-1}a_i(q^i-1) &\hspace{-6pt}=\hspace{-6pt}& q^{2m} -1
\vspace{-4pt}
\\
\displaystyle\sum_1^{2m-1}a_i &\hspace{-6pt}=\hspace{-6pt}& |\Sigma_X| \le |\Sigma_\star|=q^m+1 .
\end{array}
\vspace{-2pt}
$$

There cannot be two subspaces of  $X$   of dimension  $>m$
and $\ge m$  having zero intersection.  Thus, if    $a_k\ne0$ for some $k>m$ then
$a_k=1$ and $a_i=0$ whenever  $m\le i\le 2m-1 ,$ $ i\ne k$.  This produces the contradiction 
\vspace{1pt}
$q^{2m}-1=(q^k-1)+ \sum_1^{m-1}a_i (q^i-1) \le 
\vspace{1pt}
(q^k-1)+\sum_1^{m-1}a_i (q^{m-1}-1)
\le (q^k-1)+ (q^m +1-1) (q^{m-1}-1).
$

 \vspace{1pt}
Thus, $a_k=0$ for $k>m$, and 
\vspace{2pt}
$q^{2m}-1=  \sum_1^{m }a_i (q^i-1) \le 
 \sum_1^{m }a_i (q^m-1)
\le (q^m+1) (q^{m }-1).$
Then $a_i=0$ for $i<m$
and $a_m=q^m+1$, as required. 
\qed

\Remarks\rm
When $2m=4$ the theorem is a special case of 
\cite{CDFS,ThK} and \thmref{conics},
which suggests  the question:  \emph{Can more than one subset  like
$\Sigma_\star$   be removed in}  \thmref{using transversals}?

The last part of the proof showed that \emph{a partition of the nonzero vectors of $\F_q^{2m}$ induced by a set of proper subspaces has size at least $q^m+1,$ with equality if and only if the subspaces all have dimension $m$.}

  \section{Orthogonal  spreads}
  \label{Using orthogonal  spreads}

  Let $V$ be an  $\O^+(4m ,q )$-space (for even $q$
 and  $4m\ge8$) with quadratic form $Q$. Then $V$ has an orthogonal spread  $\Sigma$, and  $|\Sigma|=q^{2m-1}+1$ 
 (first proved in   \cite{Dillon}, then rediscovered in  
 \cite{Dye};  cf.   \cite{Ka2,KaW}).    This leads to our 
 simplest examples:
 
\Proposition
\label{orthogonal spread}
\hspace{-2pt}
$\Sigma$ is a maximal partial   spread of 
size $q^{2m-1}  +  1,$ and is  symplectic.%
   \rm

\medskip
\proof 
For even $q$, t.s. subspaces are also t.i., so    $\Sigma $ is symplectic.
{\em Maximality}:~since   $2m > 2$,
the quadratic form induced by $Q$ on any $2m$-space has a nontrivial zero.  Thus,
every $2m$-space   has  nonzero intersection with some  member of~$\Sigma $.\qed%

\Remark \rm
\label{1/2}
If    $d=2^{2 m}$ and $q=2$   then $| \Sigma|=  \frac{1}{2}d +1$.
Finding    maximal~symplectic partial spreads 
  of size $\frac{1}{2}d +1$ appears  to be  a goal of     \cite{MBGW}.   
   
\Lemma
\label{nonzero singular} 
Let 
$ E=\F_q\subseteq  F=\F_{q^k}\! $ with     $q$ even.~If 
$X$ is an $E$-subspace of an orthogonal $F$-space  and  
$|X| > q^{k^2+k},$ then $X$
contains a nonzero $F$-singular vector.

\rm
\medskip
\proof
We are given an  $F$-space  $V$ equipped with a quadratic form
    $Q$ and  associated bilinear form $f(~,~)$;
both forms have values in $F$ not in $E$.  The symbol $\perp$ will refer to the $F$-space $V$, while $\<\ ~ \ \>_L$ refers to   spanning an $L$-subspace
for $L=E$ or $F$.

For $i=1,\dots, k+1,$ we will construct 
$E$-linearly independent  
vectors   $x_1,\dots,x_{i} \in X$
and  an $E$-subspace $X_i$ 
such that
$\<x_1,\dots,x_{i}\>_E\le X_i \le \<x_1,\dots,x_{i}\>_F^\perp\cap X$    and
$|X_i|\ge|X|/|F|^{i}$.   
(In particular, $ x_1,\dots,x_{k+1}\in  \<x_1,\dots,x_{k+1}\>_F^\perp\cap X$.)

Let    $0\ne x_1\in X$ and
$X_1:=\<x_1\>_F^\perp\cap X$.  Then 
\vspace{1pt}
 $x_1\in X_1$ (since $q$ is even and hence $V$ is symplectic) and
$|X_1|=|\<x_1\>_F^\perp|| X|/| {\<x_1\>_F^\perp+X} | \ge 
|\<x_1\>_F^\perp|| X|/|V | =| X |/|F|$.

For induction,   let  $1\le i\le k$ and assume that we have $x_1,\dots,x_{i}$ and $X_i$.~Then
$|X_i| \! \ge\! |X|/|F|^i\! >\! q^{k^2+k}/(q^k)^k \ge | \<x_1,\dots,x_{i}\>_E|$.
 Let $x_{i+1}\in X_i- \<x_1,\dots,x_{i}\>_E$  and
 $X_{i+1}:=\<x_{i+1}\>_F^\perp \cap  X_i$.
Then
$x _{i+1}\in X_{i+1} \le   \<x _{i+1}\>_F^\perp\cap \<x_1,\dots,x_{i}\>_F^\perp\cap X$
and 
 $|X_{i+1}| \!  =
  \!  { |\<x_{i+1}\>_F^\perp || X_i|/| \<x_{i+1}\>_F^\perp\!  +\! X  |}\! \ge \! (|X |/|F|^i)/|F|$, as needed for induction.

Since   $\<x_1,\dots,x_{k+1}\> _E $ is  in  $ X\cap \<x_1,\dots,x_{k+1}\> _E^\perp$ 
and has size    $q^{k+1}>|F| $,
the additive map  $Q\col  
  \< x_1,\dots,x_{k+1}\>_E  \to F$   has a nonzero kernel. 
\qed  

\Remarks\rm
The preceding argument did not require anything about the nature of the quadratic form,  which  could even have a large radical.

Although the argument used the fact that all vectors are isotropic,
it can still be used for   unitary spaces and   orthogonal spaces of odd characteristic.  One minor difference is that we need to know that $X_i$ has an 
isotropic vector 
 $x_{i+1}\in X_i- \<x_1,\dots,x_{i}\>_E$.  This is clear if $X_i$
 is the span of its isotropic vectors; and that holds unless $X_i/\rad X_i$ is anisotropic, hence of dimension 1 or (in the orthogonal case)  2.  Thus, there is a 
  choice    $x_{i+1}$ for each $i$ 
 if we replace the condition $|X| > q^{k^2+k}$ by 
 $|X| > (q^2)^{k^2+k+1}$ for unitary spaces and by $|X| > q^{k^2+k+2}$ 
 for orthogonal spaces.
 
These observations do not, however, lead to useful unitary 
or odd characteristic orthogonal analogues 
of the next theorem:~unfortunately, there is no unitary spread
in dimension~$\ge6$~\cite{Thas1990}~and
no known odd characteristic orthogonal spread in dimension~$>8$.

 \Theorem
\label{nk}
If  $q$ is even   and $ m > (k+1)/2  ,$ then 
 $\F_q^{4 m k}$ has  a maximal partial spread 
of size 
$ q^{2 m k-k}+1 $
  that is  orthogonal and hence also symplectic.
\rm\medskip

\proof
Let  $V$ be an $\O^+(4 m ,q^k)$-space
with quadratic form $Q$, and let   
$\Sigma$ be an orthogonal spread in  $V$.
Let $T\col \F_{q^k}\to \F_q$ be the trace map. 
Then $Q'(v):=T(Q(v))$ is a quadratic form that 
turns $V$ into an $\O^+(4 m k,q)$-space.
Moreover,  
$\Sigma$  is still an orthogonal partial  spread in  this space,  
of size $(q^k)^{2 m-1}+1$.  

{\em Maximality}:   
If $X$ is  an $\F_q$-subspace  of $V$ of dimension $2 m k$,
then  $|X|=q^{2 m k}
>q^{k^2+k }$. By
  \lemref{nonzero singular},    $X$ contains a  nonzero 
$\F_{q^k}$-singular vector that must lie in some
member of the $\O^+(4 m ,q^k)$-spread $\Sigma$.  Thus, 
  $X$ has nonzero intersection with a member of $\Sigma$.\qed

 \question
 \label{conjectures}
Is  every 
 $\,\O^+(4 m ,q^k)$-spread  also a maximal  orthogonal partial  spread in an
  $\O^+(4 mk ,q )$-space\rm?
  \rm
This seems plausible since it is correct when either  $ m > (k+1)/2~$ 
(\thmref{nk})
 or $m=2$ \cite{Grassl} (cf. \thmref{Grassl example}(i)).
 
\Remarks \rm
If  $q=2$ and   $d=2^{2 mk}$ with $ m > (k+1)/2 $, then the    maximal symplectic   partial spreads in \thmref{nk}  
\vspace{1pt}
 have size $\frac{1}{2^k}d +1$,   resembling  \remref{1/2}.  
 Grassl's computer data \cite{Grassl} suggests much more:  
 \emph{for even $q$  there 
 appears to be  a maximal symplectic partial spread of size $2^i+1$ in $\Sp(2n,q)$-space
 whenever 
 $q^{n/2}+1\le 2^i+1\le q^{n}+1.$}

\medskip

We will need the following elementary observation several times:
\Lemma
\label{orthogonal implies symplectic}
If $\Sigma$ is a maximal  orthogonal partial  spread of  an $\O^+(4m,q)$-space
with $q$ even and 
$m\ge 2,$ then it is also a maximal  symplectic partial spread.
\rm\medskip

\proof 
Suppose that $Y\notin \Sigma $ is a t.i.~$2m$-space   such that 
$\Sigma\cup \{Y\}$ is a symplectic partial spread.
The quadratic form on $V$ restricts to a semilinear map on the t.i. subspace  $Y$;  the
  kernel is a t.s subspace  $Y_0$  of dimension 
 $\ge 2m-1$. 
 If   $  \dim \, Y_0 =2m$  then   $ Y= Y_0 $ must have the same type as the members of   $\Sigma  $   (cf.~\secref{Background}).
  
  In any case let $W$ be the t.s.~$2m$-space containing $Y_0$   having
   the same type as the members of   $\Sigma$.  By maximality,  $\Sigma\cup \{W\}$ is not an orthogonal 
  partial spread, so that $W\cap X\ne 0$ for some  $X\in \Sigma$. 
Since $\dim \, W\cap X \ \equiv 2m$ (mod~2)  we have 
$\dim \, W\cap X \ge2$.  
Since $ Y_0,$ $ W\cap X \le  W$
and $\dim \, Y_0\ge  2m-1$,  it follows that  
 $Y_0\cap (W\cap X) \ne0$ and hence that $Y\cap X \ne0.  $ \
 This contradicts the fact  that 
$\Sigma\cup \{Y\}$ is a   partial spread. \qed
 

  \section{$\O^+(4,q^k)$-space}
 \label{More maximal  symplectic partial spreads}

\Example  \rm
\label {Folklore}
  \ \em
If  $q$  is even  then an $\Sp(4,q)$-space has a maximal symplectic partial spread  of size
$ q+1 $ that is also a maximal orthogonal partial spread.\rm\medskip

\proof
An $\O^+(4,q)$-space has $(q+1)^2$ singular points partitioned by  exactly two orthogonal spreads
 $\Sigma,$ $
\Sigma^\dagger $,
 arising from the two types of t.s.~2-spaces (cf. \secref{Background});
 each member of $\Sigma $ and each member of 
 $\Sigma^\dagger $ meet nontrivially.
Possibly the most elementary (and most opaque) way to see that these are maximal symplectic spreads is to count the number of t.i. lines containing 
at least one singular point.  There are $2|\Sigma| +(q+1)^2(q-1)
=(q^2+1)(q+1)$ such 
lines, which is exactly the total  number of t.i. lines.
\qed

\Theorem
\label{Grassl example}
Let $q$ be even and $k\ge1$. 
\begin{itemize}
\item[\rm(i)]
An $\Sp(4k,q)$-space has a maximal symplectic partial   spread of size $q^k+1$   that is also a maximal orthogonal partial spread.

\item[\rm(ii)]
An $\Sp(4k,q)$-space has a maximal symplectic partial   spread of size $2q^k+1$.

\end{itemize}
\rm
 \smallskip

\proof (i)
The   preceding example  produces a 
maximal symplectic partial  spread $\Sigma$ 
of an 
$\Sp(4,q^k)$-space $V$  that is also a maximal orthogonal partial spread.
Viewed     over $\F_{q}$  
  (using a trace map as in the proof of  \thmref{nk})
   the  set  $\Sigma$ again is  
   an orthogonal partial spread.  It is a 
   maximal symplectic partial  spread by
  \cite{Grassl},
  and hence also a maximal orthogonal partial spread.      
 
 We include slightly more detail:   
in    \cite{Grassl}  the   $\F_q$-space 
   $(\F_{q^k}^2)^2$  is equipped  with the alternating bilinear form
   $\big ((u,v),(u',v') \big) :=T ( u\cdot v'-u'\cdot v)$ using the trace map 
   $T\col \F_{q^k}\to \F_q $.  The  partial spread 
   $\Sigma$
   consists of the   t.i. subspaces  
    $\{(0,0,y_1,y_2)\mid y_1,y_2\in \F_{q^k}\}$
    and
 $\{(x_1,x_2, x_2\a,x_1\a)\mid x_1,x_2\in \F_{q^k}\}$ for each
  $\a\in \F_{q^k}$.  These are  t.s.~$2k$-spaces for the  quadratic form $Q(u,v)=T(u\cdot v)$.
  In the preceding example,    $\Sigma^\dagger $  is   $ \Sigma^j$,
  where $j\col (x_1,x_2, y_1,y_2 )\mapsto 
  (x_1,y_1, x_2,y_2)$. 
  
 (ii)  
  Choose 
any $Z\in \Sigma$.  Obtain a new symplectic partial spread $ \Sigma ^\bullet$  by removing 
$Z$ and then,  for each 1-dimensional  $\F_{q^k}$-subspace $W$  of  $Z$, adjoining one 2-dimensional t.i.  $\F_{q^k}$-subspace that contains
 $W$ and is different from  both    $Z$ and   the member of  
 $\Sigma^\dagger$  containing $W$.  
 This 
  produces a maximal symplectic partial  spread of the $\F_{q^k}$-space 
  $V$   
  \cite[Remark 2.12(2)] {CDFS}.  

In fact 
{\em  $ \Sigma ^\bullet$    is also 
a maximal symplectic partial  spread of  the $\F_{q}$-space
$V$.}  For, let $X$ be a t.i.~$2k$-dimensional  $\F_{q}$-subspace of $V$ having zero intersection with all members of $\Sigma^\bullet$.  By (i),
$X$ has nonzero intersection with some member of  $\Sigma$, which therefore must be $Z$.
Then $X$ has nonzero intersection with some 
 $\F_{q^k}$-point $W$ of $Z$
 and hence with the  adjoined $\F_{q^k}$-space 
 in $\Sigma ^\bullet$
 containing  $W$, which is   a contradiction.\qed

\smallskip\smallskip
  Part    (i) contains \thmref{partial ovoid from spread}(i)
 as a special case, and
 amounts to the case $m=1$ not dealt with in   \thmref{nk}.

The proof of (i) in  \cite{Grassl} uses a neat computational idea.  
Unsupported optimism suggests that there  should also 
be a  nice geometric  proof.

\Example\rm
\label{intersect in spread}
\propref{Grassl example}(i)   points to a general construction
(compare Remarks~\ref{use ovoid}).  
Let
 $V=\F_q^{4m}$ be an orthogonal, symplectic or unitary space.
 Let $X$ and $Y$ be t.i./t.s.  $2m$-spaces with zero intersection, and let $\Sigma_X$ be a partial spread  (of $m$-spaces) of $X$.  Each $A\in \Sigma_X$ determines
  another  $m$-space $A':=A^\perp \cap Y$, and $A+A'$ is a  t.i./t.s.  $2m$-space.  Then 
  {\em $\Sigma:=\{A+ A'\mid A\in \Sigma_X\}$ is a partial spread of the same type as the underlying space $V$}.   (If $A\ne B\in \Sigma_X$ then 
  $V=X\oplus Y = (A\oplus B)\oplus    (  A'\oplus B')$,
so that   $A\oplus A'$   and $B\oplus B'$ have zero intersection.)
  
  When $\Sigma_X$ is  a  maximal partial spread (or even  a spread), some  of these partial spreads may be  maximal 
  orthogonal, symplectic or unitary partial spreads of size $q^{m}+1$
  (as in
  \propref{Grassl example}(i)
   and  \thmref{partial ovoid from spread}), but we do not see how to prove that.
  (See Question~\ref{spanning ovoid}  for instances of such symplectic partial spreads 
  that are {\em not} maximal.  
   As noted earlier, there is no unitary spread in dimension $\ge6$  \cite{Thas1990}.)%
  
\section{Projections}
\label{Projections}
Let $q$ be even.~A key ingredient of \cite{Ka2, Ka3, KaW} is the fact that there~is a natural  transition between $\O^+(4m,q)$-spreads and 
$\Sp(4m-2,q)$-spreads.  This 
uses any nonsingular point $z$ of an $\O^+(4m,q)$-space and projects into
the~symplectic space  $z^\perp/z$.  This procedure
also applies to   orthogonal and  symplectic partial spreads:

\Lemma   
\label{orthogonal to symplectic}
Let $z$ be a nonsingular point of    an $\O^+(4m,q)$-space $V$.

\begin{itemize}
\item[\rm(i)] If $\Sigma $ is a maximal orthogonal partial spread of $V, $ 
then 
$\Sigma/z:=\break
\{\<z^\perp\cap  X,z\>/z  \mid X\in \Sigma \}$ is 
a maximal  symplectic partial spread of the $\Sp(4m-2,q)$-space
 $z^\perp/z$.
\item[\rm(ii)]
If $\Sigma'$ is a maximal  symplectic partial spread of  $z^\perp/z,$
then there is a maximal orthogonal partial spread $\Sigma$ of 
$\,V\!$ such that $\Sigma'=\Sigma/z$.
Moreover$,$  $\Sigma$ is a maximal  symplectic partial spread.

\item[\rm(iii)]  If $\Sigma_1 $ is a maximal orthogonal partial spread of
 $V$  and $z_1$ is a nonsingular point of $V,$  then  
 $\Sigma/z$ and $\Sigma_1/z_1$ are equivalent symplectic partial spreads
  if and only if 
  $\Sigma_1 $ is the image of   $\Sigma $   under an automorphism of the 
orthogonal geometry of $V$ that sends    $z$  to  $z_1$.
\end{itemize}

\rm\smallskip
\proof  
(i)   \emph{$\Sigma/z$ is a symplectic partial spread}:  
If $X$ and $Y$ are distinct members of $\Sigma$ and 
$\<z^\perp\cap  X,z\>\cap \<z^\perp\cap  Y ,z\>\ne z$,  
then $z\in \<x,y\>$ for some points $x\in z^\perp\cap  X,$ 
$ y\in z^\perp\cap Y$.  Then 
 $x$ and $y$ are perpendicular to $z$ and hence to one another,
 so that $\<x,y\>$ is t.s.  whereas $z$ is nonsingular.
 
  \emph{Maximality}:~Suppose that  $(\Sigma/z ) \cup \{U/z\}$
 is a larger  symplectic partial spread for a t.i.  $2m$-space $U$
  of $V$ containing $z$.   Let $U'$ be~the hyperplane of $U$ 
  consisting~of~singular vectors
  (i.\hspace{1pt}e., $U'$ is the kernel of the semilinear map $U\to \F_q$
  induced by the quadratic form on $V$).  
 The members of $\Sigma$ all have the same type
  (cf. \secref{Background}).  
 Let~$\hat U$ be the t.s. $2m$-space of that type 
 containing $U'$.~Then $\hat U$ meets each $X\in \Sigma$ in~at most a 1-space and hence only in  0  (by
  \secref{Background},  
  $1\ge \dim\,  ( \hat U \cap X )\ \equiv 2m$ (mod~2) and hence 
$\hat U\cap  X=0$).~Thus, $\Sigma\cup \{\hat U\}$ is 
  an~orthogonal partial spread properly containing $\Sigma$, 
  whereas  $\Sigma$  is  assumed to be  maximal.
 
 \smallskip
 (ii) 
Choose a type of t.s.~$2m$-space of $V$.  If $U/z\in \Sigma'$
let $U'$ be the hyperplane of singular vectors of the t.i.~$2m$-space   $U$,
  and let  $ \hat U$ be the t.s. $2m$-space containing $U'$ of the chosen type.  Then the set 
$\Sigma$ consisting of these subspaces $ \hat U$ is an orthogonal partial spread:
distinct members of $\Sigma $  meet in at most a 1-space and hence
have  intersection 0 since  members of $\Sigma $ have the same type.  Clearly $\Sigma'=\Sigma/z$.

 \emph{Maximality}: If $\Sigma^+$ is an orthogonal partial spread
properly containing $\Sigma$, then $\Sigma^+/z$ properly contains
$\Sigma/z \! = \!\Sigma',$   
  whereas  $\Sigma'$  is  maximal.
  
The final statement follows from 
  \lemref{orthogonal implies symplectic}.

\smallskip
 (iii)  As a consequence  of Witt's Lemma \cite[p.~57]{Taylor}, an equivalence from $\Sigma /z $ to $\Sigma_1 /z _1$~lifts first to 
 $z^\perp\to z_1 ^\perp$ and then to an automorphism of the 
orthogonal geometry on $V$ 
 sending $z\to z_1$   and $\Sigma \to \Sigma_1$.  The converse is clear.\qed
  
\medskip
By (iii), a   maximal orthogonal partial spread $\Sigma$  produces many inequivalent maximal symplectic partial spreads 
 for different choices of $z$, where~the number of inequivalent ones requires knowledge of the automorphism group of~$\Sigma$.
 This was crucial in \cite{Ka2, Ka3, KaW}.
  
\Theorem  
\label{project symplectic1}
If  $k \ge2$ then there is a maximal partial 
$\Sp(4k-2,q)$-spread of size $q^{k}+1$.

\rm\medskip
\proof
Use  \lemref{orthogonal to symplectic}(i)  and
\thmref{Grassl example}(i).
\qed

\Theorem  
\label{project k}
If $ m > (k+1)/2 $ then there is a maximal partial 
$\Sp(4mk-2,q)$-spread of size $q^{2mk-k}+1$.
\rm\medskip 

\proof
Use  \lemref{orthogonal to symplectic}(i)  and
\thmref{nk}. \qed

\section{$8$-dimensional  partial spreads}
\label{$8$-dimensional  partial spreads}

In  $\O^+(8,q)$-spaces,  triality  \cite{Ti}    allows  us to use more easily visualized   points  and  partial ovoids
in place of partial spreads:
a triality map sends   orthogonal (partial) ovoids to   orthogonal (partial) spreads.
This produces  maximal  partial 
$\Sp(8,q)$-spreads when $q$ is even.

\vspace{-1pt}

\subsection{$8$-dimensional     ovoids}

Spreads and ovoids are known  in $\O^+(8,q)$-spaces 
  when  $q$  is prime, a power of 2 or~3, or 
$\equiv 2$ (mod 3) (some of these ovoids are described  in \cite{Ka1}). 
They have size $q^3+1$.  

\Lemma
\label{use orthogonal ovoid}
Let $\Omega$ be an ovoid in an $\O^+(8,q)$-space $V,$ where $q>2$.  Let $a\notin \Omega$ be a singular point   that is the only singular point in $ \<a^\perp\cap\Omega\> ^\perp$. 
{\rm (Examples  appear  below in  Appendix~\ref{appendix A}
for all even $q>2$.)}  
Then $\Omega^\bullet:= \big ( \Omega - (a^\perp\cap \Omega) \big ) \cup \{a\} $ is a 
maximal orthogonal partial ovoid of size $q^3 -q^2+1 $.%

\rm\medskip

\proof  Clearly $\Omega ^\bullet$  is  an orthogonal  partial ovoid.  If $b$ is a singular   point not perpendicular to any member of $\Omega^\bullet$ then $b^\perp\cap \Omega\subseteq a^\perp\cap \Omega$.
Since   both of these  sets have size $q^2+1$ (e.\,g., by  \cite[p.~1196]{Ka1}), we~obtain the contradiction
that  both $a$ and $b$ are  the singular point  in 
  $ \<a^\perp\cap\Omega\>^\perp$.\qed

\Theorem
\label{triality to spread}
Applying triality  $\tau $ to the preceding lemma produces a maximal orthogonal partial   spread $\Sigma^\bullet := \Omega^\bullet{}^\tau$ of size 
$q^3 -q^2+1 $
 in  an $\O^+(8,q)$-space when  $q>2$.~If $q$ is even then $\Sigma^\bullet $ is a maximal symplectic partial    spread.

In particular$,$  such maximal symplectic partial spreads exist for all even $q$.
\rm\medskip

\proof 
By the previous lemma,   $\Sigma ^\bullet$ is a
maximal  orthogonal partial       spread.  
If $q$ is even  use 
\lemref{orthogonal implies symplectic}.  \qed
\medskip
 
 We can imitate the preceding result and  remove several sets $a^\perp\cap \Omega$ by using a specific type of  ovoid. 
\theorem
\label{orthovoids}
If $q>2$ is even and $1\le s \le q/ 5,$
then an $\O^+(8,q)$-space has  maximal orthogonal partial spreads of  size 
$\,  n_s=q^3-s q^2+(s-1)(q+2)+\binom{s}{2}(q-2)+1 .
$
There is also a maximal orthogonal partial spread of size 
$n_4-1$ if $q\ge 16$.

These are also  maximal symplectic partial    spreads.

\rm\medskip
\proof  As in \thmref{triality to spread} we will construct   maximal orthogonal
ovoids.
Since this the only part of this paper involving detailed computations,   those
computations have been postponed to  Appendix~\ref{appendix A}.

For the ovoid $\Omega$ in Appendix~\ref{appendix A}, 
  \exref{examples}(i)  provides us with many  sets ${\mathcal P}$ of $s$ singular points disjoint from  $\Omega$ 
  together with  the sizes 
  $   |\bigcap_{p\in {\mathcal P} ' } p ^\perp   \cap \Omega   | $
  for  all ${\mathcal P} '\subseteq {\mathcal P}$.    Then
  \vspace{-4pt}
  $$\Omega_s^\bullet:= 
  \big (\Omega-\bigcup _{p\in {\mathcal P}}(p^\perp \cap \Omega)\big )\cup {\mathcal P}
  $$
  is an    orthogonal partial ovoid of size
  $(q^3+1)-s(q^2+1)+\binom{s}{2}(2q) -
   \vspace{1pt}
  {\binom{s}{3}(q+2)} + \cdots 
  \pm \binom{s}{s}(q+2)  +s = 
  \vspace{2pt}
  (q^3+1)-s(q^2+1)+\binom{s}{2}(2q)    -(q+2)
  + s(q+2) -\binom{s}{2}(q+2)  +(1- 1)^s(q+2)+s.  
 $ 
 \smallskip

{\em Maximality of $\Omega_s^\bullet$}:     Suppose that $b$ is a singular point 
not perpendicular to every member of
$\Omega_{s}^\bullet$.   Since   $\Omega$ is an orthogonal ovoid,
  $b^\perp \cap  \Omega$ must be contained in 
  $\bigcup_{p\in {\mathcal P}}(p^\perp \cap  \Omega ) $.  By \lemref{non-perp}, 
  \vspace{1pt}
  $s (5q- 5)\ge  \sum _{p\in {\mathcal P}} |b^\perp\cap p ^\perp \cap  \Omega | 
\ge |b^\perp \cap  \Omega|=q^2+1$,
which contradicts our assumption that $ s \le q/5$.  

The same argument can be used for    \exref{examples}(ii),
  producing the stated additional maximal orthogonal partial spreads.

Use \lemref{orthogonal implies symplectic}
 for the final assertion.  \qed
 \medskip
 
The preceding proof should be compared to   the proofs 
of  Theorem~\ref{smaller Grassl} and 
  the more elementary  Theorem~\ref{conics}.  
  In those proofs  the needed intersection sizes   are known for 
simple geometric reasons.  Here there does not seem to be a geometric explanation for the various intersection sizes   occurring in 
Appendix~\ref{appendix A}.

 \vspace{-2.5pt}

\subsection{$4$- and $5$-dimensional   orthogonal   ovoids}
The next 8-dimensional partial spreads (in \thmref
{partial ovoid from spread}) arise
from small-dimensional ovoids.  

\Example
\label {elliptic}
If $\Omega$ is an $\O^-(4,q)$-ovoid 
{\rm (i.\,e., an  elliptic quadric)}
in an $\O^-(4,q)$-space $W$ inside  a nondegenerate  orthogonal  $\F_q$-space
$V,$ then  $\Omega$   is   a maximal orthogonal partial ovoid of $V$.

\rm\medskip

\proof 
If  $x$ is any point of $V$ then $x^\perp\cap W$ contains a hyperplane of $W$ and hence contains either $p^\perp\cap W$
 for a singular point $p$ of $W$  or     $n^\perp\cap W$   a nonsingular point $n$
of $W$.  Each such hyperplane of $W$ contains a singular point of 
$W$, and hence meets $\Omega$ nontrivially.  
\qed\smallskip

A more general version of this example is a simple  consequence of  
5-dimensional
 results   \cite{Bagchi-Sastry,Ball}
 (also see \lemref{hyperplane intersections}):

\Lemma
\label {partial ovoid from ovoid}
If $\Omega$ is an ovoid in an $\O(5,q)$-subspace of  a nondegenerate  orthogonal  $\F_q$-space
$V,$ then it is   a maximal orthogonal partial ovoid of $V$.
 \rm
 \smallskip

\proof  
Once again we will show  that  
{\em each   point $x$ of $V$ is perpendicular to some 
point in $\Omega$.}  
We may assume that $U:=\<\Omega\>\not\leq x^\perp$, so that 
$H:=x^\perp \cap  U$ is a  hyperplane of $U$.  
By the preceding example, we may also assume that  
$U$ is not of type $\O^-(4,q)$.

If $H$ has type $\O^+(4,q)$ then $H$ contains a t.s. line, and each t.s. line
of $U$  meets 
each ovoid   of $U$ (by definition; see \secref{Background}).

If $H$ has type $\O^-(4,q)$ then its set $\Lambda$ of singular points is a
classical quadric.  Then 
$  \Lambda \cap    \Omega \ne \emptyset$
 \cite{Bagchi-Sastry,Ball}.

Thus, $H$ is degenerate.  If there is a singular point $y$ in its radical 
$\rad\, H$,  then every t.s. line of $U$ on $y$ meets $\Omega$ at a point perpendicular to $y$.

Finally, if  $\rad\, H $ is a    nonsingular point   then $q$ is even and the radical  
$r$ of $U$ is in $H$ (since all hyperplanes of $U$ not  containing its radical are nonsingular).    
Let ``bar'' denote the projection map $U\to U / r$.  
Then $\overline H $ is a tangent or secant plane of the  
  ovoid  $\overline \Omega$  in the 4-space  $\overline U$,  so that $\overline H$ contains 
1 or $q+1$ points of   $\overline\Omega$.  
If $T/r$  is one of these points,  then   the line $T$
has a unique singular point, and this lies in both $H\le  x^\perp$ and 
  $\Omega$.\qed

   
\question \em 
\label{spanning ovoid}
 Which ovoids in orthogonal spaces are partial ovoids in all larger-dimensional orthogonal spaces over the same field$\,?$  \rm 
This requires that all hyperplanes of the smaller   
orthogonal space
meet the ovoid.~\emph{Perhaps} this  does not   
hold for any ovoids of $\O^+(6,q)$-spaces  that span  the underlying 6-space
(and there are, indeed,  many such ovoids for which this requirement does not hold).
However, this requirement does hold  for some of the known $\O^+(8,q)$ ovoids, as   in \cite[Theorem~3.9]{Coop} 
and
for  the ovoids in \cite[Sec.~7]{Ka1} and Appendix~\ref{appendix A}.    
 
\Theorem
\label{partial ovoid from spread}
Let $q$ be a prime power.
\begin{itemize}
\item[\rm(i)] 
There are   {inequivalent}
 maximal     partial  $\O^+(8,q)$-spreads  $\Sigma$  of size 
$ q^{2 }+1 $$:$
\begin{itemize}
\item[\rm(a)] One for which   $\Sp(8,q)_\Sigma$
has a subgroup 
   $\SL(2,q^2) $~acting    $2$-transitively on $\Sigma;$ and 
   
\item[\rm(b)] One occurring when     $q$  is   odd but not prime
and for which $\Sp(8,q)_\Sigma$ is intransitive on $\Sigma.$  

\end{itemize}

\item[\rm(ii)] If $q=2^e$ then there are   {inequivalent}
 maximal     partial  $\Sp(8,q)$-spreads  $\Sigma$  of size 
$ q^{2 }+1 $
 that are orthogonal partial spreads$:$

\begin{itemize}
\item[\rm(a)] One for which   $\Sp(8,q)_\Sigma$~has a 
subgroup~$\SL(2,q^2) $~acting  $2$-transitively on $\Sigma;$ and 
   
\item[\rm(b)] One occurring when     $e>1$  is   odd 
and for which $\Sp(8,q)_\Sigma$ has a subgroup
  $\Sz(q ) $ acting  $2$-transitively on $\Sigma.$  

\end{itemize}

\end{itemize}
   \rm \medskip
   
\proof   
Let $\tau $  be a triality~map for an $\O^+(8,q)$-space $V$.
For  $\Omega$   in the preceding example or lemma,
   $\Sigma =\Omega^\tau $ is a maximal 
orthogonal partial spread of    $V$.

 For (ia) use an elliptic quadric, whose group of isometries produces the
 last part.  For (ib) there are other choices for $\Omega$ in 
 \lemref{partial ovoid from ovoid}, such as those in \cite[Sec.~5]{Ka1}.

If $q$ is even then \lemref{partial ovoid from ovoid} applies, where the  only known choices for   $\Omega$   
 are an elliptic quadric (\exref {elliptic})
and a Suzuki-Tits ovoid (see 
Appendix~\ref{Background on Suzuki-Tits  ovoids}).  
The stated groups arise from subgroups  of  $\Omega^+(8,q)$
 acting on $\Omega $.  
 
The various partial spreads  are inequivalent as orthogonal partial spreads, 
since the corresponding maximal   orthogonal partial ovoids $\Sigma ^{\tau^{-1}}=\Omega $ 
 are inequivalent.
 However, this is orthogonal inequivalence, which is not the same as symplectic inequivalence in (ii).  
 
 Nevertheless, the symplectic partial spreads in (iia) and (iib) are inequivalent. 
 This can be proved using   the groups appearing in (iia) and (iib),
 but a geometric proof  is simpler.   The construction of the  orthogonal partial spread in (iia) provides us with   
 $q^2+1$ t.i.~4-spaces (in fact, t.s.~4-spaces in 
 $\Omega^\perp{}^\tau$) meeting each of its members.  If the symplectic partial spreads are equivalent then 
 (iib) has $q^2 +1$    t.i. 4-spaces  $U$  meeting each of its members.  Then $U$ is t.s. since it is spanned by singular vectors.
 Now  $(\Sigma\cup \{U\})^{\tau^{-1}} =\Omega\cup  \{U^{\tau^{-1}}\} $  for 
 a singular point $ U ^{\tau^{-1}}$ 
 of $ \Omega ^\perp$.  This contradicts the fact that, for  (iib),   
 the 3-space
$ \Omega ^\perp$ contains only $q+1$ singular  points.
 \qed
  
\Remarks\rm
\label{use ovoid}
We   excluded $q=2 $  in (iib)~since that produces the same partial spread as
in~(iia).
Part (iia) is a very special case of a result in   \cite{Grassl} (cf. \thmref
{Grassl example}).  Is there an analogous generalization of (ii)?

Note that $\Sp(8,q)_\Sigma$ even contains subgroups 
 $\SL(2,q^2) \times   \SL(2,q^2) $ in (i) and  $\Sz(q) \times \O(3,q) \cong
 \Sz(q) \times \SL(2,q) $ in (iib).

In both (i) and (ii)  there are t.s.  4-spaces $X,Y$ such that 
 the members of 
$ \Sigma$ meet $X$ and $Y$ in   spreads of each (cf. 
\exref{intersect in spread}).

See \cite{PW} for a survey  of  
$\O(5,q)$-ovoids.

\subsection{Extending a partial spread}
\label{Extending a partial spread}
\emph{How can one search for maximal symplectic   partial spreads}?
One obvious answer is to start with a symplectic or orthogonal partial spread and try to extend it to a maximal one
(this was the computational method used to produce the table in \cite{Grassl}).   The
instances  considered below may 
have  extensions to maximal ones other than the ones we provide.

Once again, points are easier to deal with than subspaces.

\subsubsection{$\O^-(4,q)$-ovoids}
\label{$O-(4,q)$-ovoids}
 A  simple example of an
  orthogonal partial ovoid is $(\Omega -\{p\})\cup\{x\}$, where 
  $p$ is a point in
the set   $\Omega$    of singular points of
an $\O^-(4,q)$-space~$U$   and  $x\notin U$ is a singular point 
in   $ (p^\perp\cap U)^\perp -U^\perp$.

 \Lemma
 \label{extend to maximal partial ovoid}
 For any $q$ an $\O^+(8,q)$-space has a maximal orthogonal partial ovoid of size $2q^2+1$.
   \rm 
   \medskip

\proof In an $\O^+(8,q)$-space $V$ consider anisotropic 2-spaces $A,A'$ and a totally singular 2-space $\<p,p'\>$ such that
$\<A,A',p,p'\>= A\perp A' \perp  \<p,p'\>$.  
Let 
$E= \<A, p \> $ and $ E '= \<A', p' \>, $ and let   $ x$ be a point 
    \vspace{1pt}
    of $ \<p,p'\>-\{p,p'\}$. 
    
    Let $U $ and $U' $ be non-perpendicular $\O^-(4,q)$-subspaces 
of $V$ such that
$E'^ \perp >  U>  E $ and 
$E^ \perp > U'>  E'$. 
(In order to construct these, note that  $p$ and $p'$ are in t.s. lines $\ne \<p,p'\> $ of the $\O^+(4,q)$-space   $ ( A\perp A')^\perp$.  Choose singular points 
$u,u'  \in     ( A\perp A')^\perp - \<p,p'\>$  perpendicular to $p'$ and $p$, respectively, but not to each other.  
Then $U:= A\perp \<p,u\>=\<E ,u\>$ and
$U':= A'\perp \<p',u'\>=\<E' ,u' \>$  are non-perpendicular $\O^-(4,q)$-subspaces such that 
$U=\<A,p,u  \><\<A' ,p ' \>^\perp=E'{} ^\perp$ and 
$U' < E  ^\perp$ behave as required.)

If $\Omega$ and~$\Omega'$ are the sets of singular points of $U$
 and $U '$, respectively, we claim~that 
 $$
 \Omega^\bullet:= ( \Omega -\{p\}   )\cup 
  ( \Omega' -\{p' \}   )\cup \{x\}
 $$
 behaves as stated in the lemma.   Clearly,  $| \Omega^\bullet|=q^2+q^2+1$.
 
 \smallskip
 \emph{Orthogonal partial ovoid}:
 $x^\perp \cap U=p^\perp \cap U=E$ has only one singular point $p$, and $p\notin 
 \Omega^\bullet$.  Suppose that there are perpendicular
 singular points  $y\in \Omega-E $  and   $y'\in \Omega ' -E ' $.
 Since $y\in U < E'^\perp $ and $y' <E^\perp $, while $E$ and $E'$ are perpendicular, we obtain  to the contradiction that 
 $\<y,E\>=U$ and $\<y',E'\>=U'$ are perpendicular.%
 \smallskip
 
 \emph{Maximality}:
 Suppose that $h$ is a singular point such that 
 $ h^\perp \cap \Omega^\bullet =\emptyset$.   Then  
 $h^\perp\cap U$ is a hyperplane of $U$ and hence contains a singular point, which must be $p$.
Then   $h^\perp\cap U = p^\perp\cap U =E $.~Also 
 $h^\perp\cap U '= E'$.  Now $h\in \<E,E'\>^\perp = \<p,p'\>$,
 which contradicts the assumption 
 that $h$ is not perpendicular to $x\in \Omega^\bullet$. \qed
  
 \Theorem
 \label{2q2+1}
 For any $q$ an $\O^+(8,q)$-space has a maximal orthogonal partial spread $\Sigma$ of size $2q^2+1$.   If $q$ is even then  $\Sigma$ 
 is symplectic.
   \rm 
   \medskip

\proof 
Applying triality to  the lemma proves the first part, while 
\lemref{orthogonal implies symplectic} implies the second
part.  \qed
\medskip

When $q$ is even,    \thmref{Grassl example}(ii) contains  another maximal symplectic partial spread of size $2q^2+1$  that
  need not be orthogonal.
  
Note that these examples, and others earlier in this section, 
would have been awkward to describe using t.s. 4-spaces instead of points.

\subsubsection{Suzuki-Tits~ovoids} 
\label{Using a Suzuki-Tits ovoid}

Another example of an
  orthogonal partial ovoid is $(\Omega -\{p\})\cup\{x \}$, where 
  $p$ is a point of   a Suzuki-Tits ovoid    $\Omega$   in 
an $\O(5,q)$-space $U$   and  $x \notin U$ is a singular point 
in   $ (p^\perp\cap U)^\perp -U^\perp$ (see 
Appendix~\ref{Background on Suzuki-Tits  ovoids}).  This time it is easier to extend this to  a maximal orthogonal partial ovoid
of an $O^+(8,q)$-space.
In the next section we will   see further advantages of  $\Omega$
   over an elliptic quadric.

 \Theorem
 \label{q2+q+1}
 If  $q=2^{2e+1}>2$ then  an $\O^+(8,q)$-space has a maximal orthogonal partial spread $\Sigma$ of size $q^2+q+1$
 that is symplectic.
   \rm 
   \medskip

\proof  By    triality and    \lemref{orthogonal implies symplectic},
we need to construct a maximal orthogonal partial ovoid of the stated size
in an   $\O^+(8,q)$-space $V$   containing $U$.
The radical $r$   of $U$  is also the radical of the 3-space   $U^\perp $,
and $(p^\perp\cap U)^\perp = \<p,U^\perp \> = p\perp U^\perp $
for $p\in  \Omega$.
Each singular point in the 4-space $(p^\perp\cap U)^\perp$ lies on a t.s. line
containing $p$ and meeting $U^\perp $ in one of its $q+1$ singular
 points.
 
 For each singular point $x _0$ in $ U^\perp $ let $x$ be any  point
 in    $ \<p,x_0\>  -   \{ p,x_0\}$.  Let $X$ be the resulting set of   
 $q+1$ points $x$.   We claim~that \vspace{-4pt}
 $$
 \Omega^\bullet:= ( \Omega -\{p\}   )\cup  X\vspace{-4pt}
 $$
 behaves as required.   Clearly,  $| \Omega^\bullet|=q^2+q+1$.
 
 \smallskip
 \emph{Orthogonal partial ovoid}:
  $x^\perp \cap U=p^\perp \cap U $ 
 since $x_0^\perp\ge U$.   Then~$x^\perp \cap\Omega=\{p\} .$  
 No two members of $X$ are perpendicular since no two singular points in  $U^\perp$ are.
 \smallskip  
  
 \emph{Maximality}:
 Suppose that $h$ is a singular point such that 
 $ h^\perp \cap \Omega^\bullet =\emptyset$. 
 Then  
 $h^\perp\cap U$ is a hyperplane of $U$ that cannot contain a point of
$\Omega -\{p\}$. 
By   \lemref{hyperplane intersections}, 
 $h^\perp\ge h^\perp\cap U = p^\perp\cap U$.
 Then $h\in (p^\perp\cap U) ^\perp = \<p,U^\perp \>$, 
 so  that $h$ lies on one of the above lines $\<p,x_0\>$,
 whereas  $ h^\perp \cap X =\emptyset$.  \qed%

 \questions   Instead of using a single pair $(p,x)$ for replacement what happens if several such pairs are used\emph{?}  Can {\rm\secref{$O-(4,q)$-ovoids}} be handled better than at present in order to use several replacement pairs\emph{?} 
 
\subsection{Small maximal partial spreads }
\label{Grassl's Conjecture}
We will describe counterexamples to Grassl's conjecture, which was stated in the Introduction.
  Grassl \cite{Grassl} has  also found counterexamples to his conjecture 
    in an $\Sp(8,8)$-space
    by a  computer search.

\Theorem
\label{easiest Grassl}
If $q=2^{2e+1}>2$ then there is a maximal   partial 
$\O^+(8,q)$-spread of size $q^2-q+1;$
this is also a maximal  partial $\Sp(8,q)$-spread.
\rm\medskip  

\proof
In view of   triality  and  \lemref{orthogonal implies symplectic}, 
it suffices  to construct a maximal   partial 
$\O^+(8,q)$-ovoid of size $q^2-q+1.$
We use  the  notation in \secref{Using a Suzuki-Tits ovoid} and Appendix~\ref{Background on Suzuki-Tits  ovoids}.

Let $\Omega$ be a Suzuki-Tits ovoid in an $\O(5,q)$-space $U$.
Embed $U$ into an $\O^+(8,q)$-space $V$.  

 Let
   $\Omega^\bullet:=
{\big (\Omega-(x^\perp \cap \Omega)  \big )\cup\{x\}}$
for a singular point $x$  of $U$ not in $\Omega$ (this uses 
 $\dim \,U>4$).  Then  $|\Omega^\bullet|=q^2-q+1$ 
 and $ \Omega^\bullet $  is an orthogonal partial ovoid of $U$  and hence   of $V$.

 \smallskip  
  
 \emph{Maximality}:
Suppose that $h$ is a singular point of $V$ such that 
$h^\perp\cap\Omega^\bullet =\emptyset$.  
We will consider the possibilities for the  hyperplane  $h^\perp \cap U$  of $U$ in  \lemref{hyperplane intersections}.   
We have  $  h^\perp \cap\Omega 
\subseteq x^\perp \cap \Omega $
since $h^\perp\cap\Omega^\bullet =\emptyset$. 
Also,   $\Omega^\perp =U^\perp < x^\perp$
since $x \in U=\< \Omega\>$.%
\smallskip

Case 1. $h^\perp  \cap \Omega =\{p\}$ for some 
$p\in x^\perp \cap \Omega.$ 
Then    
$ h^\perp \cap U =  p^\perp \cap U $
since  \lemref{hyperplane intersections} implies that  $p^\perp \cap U $ is the only hyperplane of $U$ meeting $\Omega$   just in $p$.  
Then $h^\perp\ge h^\perp \cap U =  p^\perp \cap U $, so that
$h\in \<p  , U^\perp\> \le x^\perp$, 
whereas  $h$ is assumed not to be perpendicular to $x\in \Omega^\bullet$.

\smallskip
Case 2. 
$|h^\perp   \cap \Omega|=q+1.$   Since 
$h^\perp \cap\Omega \subseteq x^\perp \cap \Omega $
for  sets of size $q+1$,
we have    $h^\perp\ge \<h^\perp \cap\Omega\>= 
\<x^\perp \cap \Omega\> $,
where  $\<x^\perp \cap \Omega\>= 
x^\perp \cap U $ by the end of \lemref{hyperplane intersections}(ii).
  Then 
$h\in  \<x  , U^\perp\> \le x^\perp$, which
   produces   the same contradiction   as before. 
(This is where  an elliptic quadric  $\Omega$   would not suffice:
  we would only have    $\<x^\perp \cap \Omega\> <x^\perp \cap U $
since $\dim \, U=5$.)

\smallskip
Case 3. $1< |h^\perp  \cap \Omega |  <q+1$.  
Since   $  h^\perp  \cap \Omega$ lies in a set $x^\perp \cap \Omega $ that projects into a plane of 
 $U/ r$,  this contradicts the irreducibility in    \lemref{hyperplane intersections}(iii). 
\qed

\medskip
We can go further  (mimicking  the proofs of 
Theorems~\ref{orthovoids}
and \ref{conics}):

\Theorem
\label{smaller Grassl}
An $\O^+(8,q)$-space has a maximal orthogonal partial spread of size 
$q^2-sq+2s-1$  whenever $q=2^{2e+1}>2$ and 
$1 <  s\le \sqrt{q/2} -1 $.  Each of these is a maximal partial symplectic spread.

In particular$,$ there is a maximal partial $\Sp(8,q)$-spread
of size $q^2-\sqrt{q^3/2} +q+\sqrt{2q}-3.$ 
  \rm  \smallskip

\proof 
Once again we will construct maximal partial
$\O^+(8,q)$-ovoids. Let  $\Omega$, $U$, $r$   and $V$ be as before.   
Choose  distinct $a,b \in \Omega$. 
Then  $\{a,b\}^\perp \cap U$ is a nondegenerate 
plane containing   $r$.  There are $q+1$ singular points   
$x\in \{a,b\}^\perp $.  These produce $q+1$ subspaces 
$\<x,a,b,r\> =x^\perp \cap U$ that induce a partition
 of $\Omega-\{a,b\}$ 
using the $q+1$    {circles}   $\Omega_x:=x^\perp \cap \Omega$ of the inversive plane $\I(\Omega)$  determined by $\Omega$   \cite[Sec.~6.4]{Demb}.

Let $\mathcal X$ be any set of $s $ singular points 
$x\in \{a,b\}^\perp $.
  We will show  that   
  \emph{
  $$\Omega ^\bullet :=
   \big (\Omega- \bigcup_{x\in\mathcal X }  \Omega_x \big )\cup  \mathcal X
    \vspace{-4pt}
   $$
 is a maximal partial ovoid of the stated size.  }

\smallskip
1.  $|\Omega ^\bullet |= 
(q^2+1) -2 - |\mathcal  X| (q-1) +  |\mathcal  X|  $.

\smallskip

2.  {\em Partial ovoid}:  If $x\in \mathcal  X$ then     
$\Omega_x= x   ^\perp\cap \Omega $
was replaced by   $ x  $,  and $\mathcal  X$ lies
   in a conic of  $\{a,b\}^\perp \cap U$.
   
\smallskip

3.  {\em Maximality}: 
  Suppose that $h$ is a singular point of $V$ such that 
  $h^\perp\cap \Omega^\bullet=\emptyset.$
Then $h^\perp \cap\Omega 
\subseteq   \cup_{x\in \mathcal  X}\Omega_x $.
We will consider the various possibilities in 
\lemref{hyperplane intersections}
for the hyperplane  $h^\perp \cap U$   of $U$.
\smallskip

Case 1.   $h^\perp  \cap \Omega =\{p\}$ for some 
$p\in \Omega.$  Then $p\in \Omega_x $ for some  
$x\in \mathcal X \subset   \Omega^\bullet$. 
By  \lemref{hyperplane intersections},   
$h^\perp \cap U =  p^\perp \cap U $.  This is inside  $h^\perp$, so that
  $h\in \<p  , U^\perp\> \le x^\perp$,   
whereas  $h$ is assumed not to be perpendicular to $x\in \Omega^\bullet$.

\smallskip
Case 2.
$h^\perp  \cap \Omega$  is a circle. If $h^\perp  \cap \Omega$ contains
$\{a,b\}$ then $h^\perp  \cap \Omega= \Omega_x\subset x^\perp$
for some $x\in \mathcal X$ since the circles 
$ \Omega_y,$  $ y\in  \{a,b\}^\perp$, induce a partition of 
$\Omega-\{a,b\}$.
Then $h^\perp $ contains $\<\Omega_x\cap  U\> =x^\perp  \cap  U$ 
by \lemref{hyperplane intersections}(ii), 
 which   again produces the  contradiction  $h\in
 \<x,U^\perp\> \le  x^\perp$. 
 
 If $h^\perp  \cap \Omega$ does not contain $\{a,b\}$ then it meets
 each circle $\Omega_x$, $x\in \mathcal X$, in at most two points.  This produces the contradiction
 $q+1=|h^\perp  \cap \Omega|\le 2| \mathcal X |=2s$.

\smallskip
Case 3. 
$h^\perp  \cap \Omega$  is an  \ orbit of a cyclic group $T<G$ of order $|h^\perp  \cap \Omega|=q\pm \sqrt{2q}+1$
(\lemref{hyperplane intersections}(iii)).  We  replace the argument used in \thmref{easiest Grassl} by  counting helped by $T$.
Note that $ |T |$ divides $ q^2+1$  and hence is relatively prime to $q(q-1)$,  the  order of
the stabilizer in $G$ of a circle \cite[Theorem~9]{Suz}.  
Thus,
given circles 
$C_1$ and $C_2$, at most one   element  of $T$ can  send  $C_1$ to  $C_2$.  

For each $t\in T$ we have 
$h^\perp  \cap \Omega=(h^\perp  \cap \Omega)^t \subseteq \cup_{x\in \mathcal X}\Omega_x^t$,
involving  two  sets of $s$ circles:
$ \{\Omega_x  \mid x\in \mathcal X\}$ and 
$ \{\Omega_x^t \mid x\in \mathcal X \}$.  For an ordered pair $x,y$
of distinct elements of  $\mathcal X$ there is at most one 
such  $t\ne 1$ with $\Omega_x^t=\Omega_y$.
  Thus, if we choose $t$  to be one of 
at least  
$|T|-1-s(s-1) \ge q - \sqrt{2q}-s(s-1) >0$ elements of $T$ that 
do not behave this way for all $x,y$,  then   we will have two disjoint sets of $s$ circles, with  the union of each set containing $h^\perp  \cap \Omega$.  Since distinct circles meet in at most two points,  
$q\pm \sqrt{2q}+1 =| h^\perp  \cap \Omega | \le s\cdot s\cdot2$, which is not the case.

\smallskip
Case 4.~$\Lambda:= h^\perp  \cap \Omega$ has size $q+1$,
and its  stabilizer in $G$  
has a  cyclic subgroup $T$ of order $q-1$  having orbit-lengths $1,1,q-1$ on 
$\Lambda$ (\lemref{hyperplane intersections}(iv)).
There are $q+1$  orbits of $T$ of size $q-1$; 
every nontrivial element of $T$ fixes just two  points.

 If $a$  and  $b$ are not the two points fixed by $T$ then it follows that 
every nontrivial element of $T$ moves every $\Omega_x$.
Thus,  the argument in Case 3 can be repeated
(this time with $|T|-1-s(s-1) = q-s(s-1) >0$ and the contradiction
$q +1  \le s\cdot s\cdot2$).

If $T$ fixes $a$ and $b$ then $\Lambda$ is not one of  the circles
$\Omega_y$ by  \lemref{hyperplane intersections}(iv).  
Note that $T$ fixes two circles containing $  \{a,b\}$  and is transitive  
on the remaining $q -1 $ circles containing $  \{a,b\}$. 
Thus, if a nontrivial element  of $T$   fixes   
$\Omega_y$  then $\Omega_y-\{a,b\}$ is an orbit 
 of $T$,
while $\Lambda-\{a,b\}$ is a different orbit, so that $\Omega_y$ can be deleted in our 
union (of $\Omega_x$, $x\in \mathcal X$)  that contains $\Lambda$.
Since  $T$ is transitive on the set of $q-1$
circles $\Omega_x$ that it does not fix, we obtain the contradiction
 $s=|\mathcal X|\ge |T|=q-1$.
(Alternatively,   the argument in Case 3 can be repeated again.)
 \qed
\smallskip

We have   proved, more generally, that 
\emph{$\Omega^\bullet$ is a maximal partial ovoid of any nonsingular orthogonal  $\F_q$-space  containing $U$} since every hyperplane of 
$U$ has nonempty intersection with $\Omega^\bullet$
(cf. Question~\ref{spanning ovoid}).

\subsection{$\Sp(6,q)$-space consequences}
\label{$Sp(6,q)$-space consequences}

 \Theorem
 \label{Sp6 corollary}
 For even $q>2, $ an $\Sp(6,q)$-space has maximal symplectic partial spreads of  size
 \begin{itemize}
  \item[\rm (i)] $n_1 = q^3-q^2+1 ,$ 
 \item[\rm (ii)] $q^2+1,$  
 \item[\rm  (iii)]  $2q^2+1,$
  \item[\rm  (iv)]  $q^2+q+1$ if $q=2^{2e+1},$
  
    \item[\rm  (v)]  $q^2-q+1$ if $q=2^{2e+1},$

  \item[\rm  (vi)]      $q^2-sq+2s-1$  if $q=2^{2e+1} $ and 
$1\le s\le \sqrt{q/2} -1 ,$

  \item[\rm  (vii)]  
 $n_r$ if $1\le r\le q/ 5 $  
{\rm(where $n_r$ is defined in \thmref{orthovoids}),}
and
 \item[\rm  (vii)]  $n_4-1$ if $q\ge16$.
\end{itemize} 
\medskip\rm

\proof
Use 
  \lemref{orthogonal to symplectic}(i) together with  
  Theorems~\ref{triality to spread},
  \ref{partial ovoid from spread},
    \ref{2q2+1},  
     \ref{q2+q+1},
 \ref{easiest Grassl},
  \ref{smaller Grassl}
 and   \ref {orthovoids}.  \qed
  
\Example\rm
By \lemref{orthogonal to symplectic}(ii), 
the set of sizes of maximal  partial $\Sp(6,4)$-spreads  
is contained in 
the set of sizes of maximal  partial $\Sp(8,4)$-spreads.
This can be compared with the list in  \cite{Grassl}.

\section{$6$-dimensional  partial spreads}
\label{Using groups}

We   again consider   arbitrary characteristic.
In characteristic 2 the  examples   in the next theorem   already
appear  in 
\thmref{Sp6 corollary}(i)  but using an entirely different method to prove maximality. 
 
\theorem
\label{group}
If  $q$ is a prime power  then
an $\Sp(6,q)$-space  has a maximal symplectic partial spread  of  size 
$q^3-q^{2}+1$. 

\rm\medskip
\proof      
In an $\Sp(6,q)$-space let $\Sigma$ be a desarguesian spread 
   preserved by   $G= \SL(2,q^3)=\Sp (2,q^3) <\Sp(6,q)$.
   Let $X\in \Sigma$.
 Let $U$ be a t.i. 3-space such that $U\cap X=L$ is a line.  
 If  $\Sigma_U$ is the set of members of 
 $\Sigma$  met nontrivially by $U$,
 then  we will show that     \emph{$\Sigma^\bullet := (\Sigma -\Sigma_U)\cup\{U\}$
 is a maximal symplectic partial spread
 of   size $q^3-q^2+1$. }
   
 If $U$ meets $Y\in \Sigma-\{X\}$ nontrivially then $U\cap Y$ must meet
 $U\cap X=L$ trivially and hence is a point; the number of such points is
 the number   $q^{2}$  of points in $U$ not in $L$. 
 Thus,   $|\Sigma_U |= q^2+1$ and    $\Sigma^\bullet $
 is a symplectic partial spread
 of   size $q^3-q^2+1$.  
 
The set-stabilizer $G_X$ of $X$
 has order $q^3(q^3-1)$, with  an abelian normal subgroup 
 $Q$ of order $q^3$  inducing $1$ on $X $
 and  a cyclic subgroup  $S$ of order $q^3-1$ that is transitive on  both $X-\{0\}$  and  the $q^2+q+1$ lines $L$ of $X$. 
     Then   $|G_L|=q^3(q-1)$.
Since   $Q$  is  transitive
 on the $q$ t.i.~3-spaces $\ne X$ containing $L$, 
     $|G_ U |=q^2(q-1)$ 
      and  $Q_U$ fixes each of those  3-spaces.%

  Since   $S$  is   transitive on both  
the $q^2+q+1$ lines $L$ of $X$ the 
   $G_X$-conjugates  of $Q_U$, we obtain a bijection $L\mapsto Q_U$ 
   between these sets.

Also   $Q$ is transitive on 
 the $q^3$ points in $\{ L^\perp\cap  Y \mid Y\in \Sigma-\{X\} \}$.
 Since~$U$ contains $q^{2}$ of these points, and each such point and 
 $L$ generate $U$, it follows that~$Q_U$ is transitive on these
 $q^{2}$  points.  
  Since  $Q_U$  fixes each t.i.~3-space   containing 
 $L$,   the  $q$  orbits of  $Q_U$ partitioning  
 $ \Sigma-\{X\}$ correspond  to  the $q$
  t.i.~3-spaces~$\ne X$  containing~$L$.%
  
     \emph{Maximality}:  Assume that $W\notin \Sigma^\bullet$
 is a t.i.~3-space such that $\Sigma^\bullet \cup\{W\}$ is a symplectic partial spread.  
 Then   $\Sigma_W\subseteq \Sigma_U$  since $\Sigma$ is a spread.
 Clearly, 
   $W$ meets each member of $\Sigma$ in 0, a point or a line.
 Since $|\Sigma_W|\le |\Sigma_U|={q^2+1}$, some intersection~is a line, and
  it is unique (since two lines of $W$ would meet
 nontrivially).~Thus, $W$ arises in the same manner as $U$,
 and $G_W$  acts on ${\Sigma_W=\Sigma_U}$.
 
 We cannot have $   L=X\cap U = X \cap W  $ in view of the above orbit 
 partition, and we cannot have 
   $  G_U = G_W  $ in view of the above bijection.
 Then $\< G_U,G_W \>$
  is generated by   distinct subgroups of 
 $G$ of order $q^2(q-1)$, and hence has a subgroup of order $q^3$
 with an orbit on $\Sigma$ of size $q^3$.
 (This is  a very special case of 115-year-old group theory
 summarized in  \cite[Ch.~XII]{Dic}.)
 This contradicts the fact that $|\Sigma_U|<q^3$.~\qed

\section{$4$-dimensional  partial spreads}
\label{Generalized quadrangles}

Finally, we survey     families   of maximal partial spreads of $\Sp(4,q)$-spaces.
See \cite{CDFS,Grassl} for lists and tables  of  known families.   
As suggested in  
   \secref{Introduction}, we can use more easily visualized points in $\O(5,q)$-space
instead of lines in $\Sp(4,q)$-space due to the Klein correspondence 
\cite[p.~196]{Taylor}.
 
\Theorem [{\cite[p.~1940]{CDFS}},   
{\cite[Theorem~6.6]{ThK}}]
\label{conics}
\        \vspace{7pt}     
\begin{itemize}
\item[\rm (i)]
\vspace{-8pt}  
For odd $q$  an  $\Sp(4,q)$-space has  a maximal partial spread of size
$q^2-s q+ 3 s-1 $  whenever  $1\le  s<(q + 1)/2$.

\item[\rm (ii)]
For even $q$  an  $\Sp(4,q)$-space has  a maximal partial spread of size
$q^2- sq+ 2 s -1 $  whenever  $1\le  s<(q + 1)/2$.
\end{itemize}
\rm\smallskip

\proof We will construct maximal partial $\O(5,q)$-ovoids.  Start with an $\O^-(4,q)$ 
ovoid $\Omega$ in a 4-dimensional subspace $U$.
Choose  distinct $a,b \in \Omega$. 
    If $y$ is a singular  point not in $\Omega,$ then $y^\perp\cap \Omega$ 
is an oval (more precisely,  a conic in $\<y^\perp\cap \Omega\>$).%
\smallskip

(i)  
The planes $E$ of $U$ containing $\{a,b\}$ fall into two sets 
$\Pi_k$, $k=0$ or 2, each of size  $(q+1)/2$,
 such that the nonsingular line 
$E^\perp $ has exactly $k$  singular points.
Let $\Omega_E:=E\cap \Omega $.
 With each  $E\in \Pi_2$   
 is an associated $E'\in \Pi_2$ (for an involution $E\mapsto E'$
 without fixed points) such that
there  are exactly two singular points 
 $x_{E },x_{E' } $ in 
$ E\cap \{a,b\}^\perp $  and
$(*)\ \Omega_E = 
x_{E } ^\perp\cap \Omega =
 x_{E '} ^\perp\cap \Omega  =\Omega _{E '}$.
Then $\{x_{E },x_{E' }  \mid E\in \Pi_2 \} $ lies in  a conic in the  plane
$\{a,b\}^\perp $.  
The members of 
$\Pi_  0\cup \Pi_ 2$ induce a partition of $\Omega-\{a,b\}$.  
 
Let $\mathcal S$ be any set of $s < (q+1)/2$ planes $E \in \Pi_ 2 $ such that the conics~$\Omega_E$,$\,E\in\mathcal S$, are distinct
(cf.~$(*)$). 
  We claim that   
  \emph{$\Omega ^\bullet :=
 { \big (\Omega- \bigcup_{E\in\mathcal  S }  \Omega_E \big )\cup \bigcup_{E\in\mathcal S}\{x_{E},x_{E' }\} }$
 is a maximal  partial $\O(5,q)$-ovoid of the stated size.  }
 
 \vspace{2pt}

1.  $|\Omega ^\bullet |= 
(q^2+1) -2 - |\mathcal  S| (q-1) +2 |\mathcal  S|  $.

2.  {\em  Orthogonal partial ovoid}:  If $E\in \mathcal  S$ then     
$\Omega_E= x_{E}  ^\perp\cap \Omega 
= x_{E' } ^\perp\cap \Omega $
was replaced by   $\{x_{E},x_{E' }\}$, 
lying  in a conic of  $\{a,b\}^\perp$.

3.  {\em Maximality}:  Every point of $\Omega$ is either in  
$\Omega ^\bullet$ or is perpendicular to $\{x_{E},x_{E' }\}$
for some  $ E \in\mathcal S$.
 Suppose that $h\notin \Omega$ is   a singular point such  that
 $h^\perp\cap  \Omega ^\bullet =\emptyset$.
Then $h^\perp\cap \Omega\subseteq  \bigcup_{E\in\mathcal S}  \Omega_E$
and $h^\perp\cap \Omega$ is either a point or a circle of 
the inversive plane   
$\I(\Omega)$  determined by $\Omega$   \cite[Sec.~6.4]{Demb}
(compare  \exref{elliptic}).

If  $h^\perp\cap \Omega$ is a point $p$ then
$h^\perp\cap U$ is the tangent plane to $\Omega$ at $p$  in $U$.
Then 
$h\in (h^\perp\cap U)^\perp = (p^\perp\cap U)^\perp =\<p,U^\perp\>$,
which has just one singular point $p$, whereas $h\notin \Omega$.  Thus, 
$h^\perp\cap \Omega \subseteq  \bigcup_{E\in\mathcal S}  \Omega_E$
 is a circle.  
 
If $h^\perp\cap \Omega=\Omega_E$ with $E\in\mathcal S$,  then 
$h\in (\Omega_E)^\perp =  E^\perp =  \<x_{E},x_{E' }\>$, whereas 
$h\notin \{x_{E},x_{E' }\}$.
Thus,   $h^\perp\cap \Omega$ is a circle lying  in the union of 
$s$ other circles, each of  which it meets at most twice.  This produces 
the contradiction $q+1=| h^\perp\cap \Omega|\le 2s$.%

 \smallskip
(ii)  This is proved as above but  is simpler:  
${(E^\perp\cap \Omega)^\perp}$ contains just one singular point
for each plane $E$ of $U$ containing $\{a,b\}$;  no 
permutation $E\mapsto E'$ is involved.\qed

 \medskip
    
  It is tempting  to hope that the above argument only used an $\O(5,q)$-ovoid $\Omega$.
 However, when $q$ is odd the intersections $x^\perp \cap \Omega$
 are not sufficiently well-behaved.  When $q$ is even       \secref{Grassl's Conjecture} contains
   versions of the preceding argument in $\O^+(8,q)$-space
  (also compare Theorem~\ref{orthovoids}).
  
\Example  \rm
\label{3q-1}
 A maximal partial ovoid of size $3q-1$
in    $\Sp(4,q)$-space, $q\ge4,$ is constructed in 
\cite[p.~1939] {CDFS}.  The proof in that paper shows that this is 
a maximal partial ovoid 
in   $\Sp(2m,q)$-space for all $m\ge2$.  

This partial ovoid is the set of  points in 
$\bigcup_1^3\big (\<x_i,y_{i+1}\>-\{x_i,y_{i+1}\}\big)\cup \{x,y\}$ (subscripts mod~3), 
where $x_1, x_2, x_3, x$ are four points of   $X$
and 
$y_1, y_2, y _3, y$  are four points of  $Y$ for 
 t.i.~2-spaces $X,Y$ intersecting in 0, with  each
 pair  $x_i,y_i$   perpendicular and $x,y$   not  perpendicular.

Dualizing \cite[p.~196]{Taylor}  produces a maximal symplectic partial spread 
of size $3q-1$
in    $\Sp(4,q)$-space for even $q\ge4$.  

\Example  \rm
 For arbitrary $q$,
  \cite{RS,PRS}  has 
  integer intervals that consist of sizes of maximal
   partial  $\Sp(4,q)$-spreads. 
 
 \Example  \rm
 There is a maximal partial spread  of size $q^2-1$
in    $\Sp(4,q)$-space  for $q\in\{3,5,7 ,11 \}.$   This is constructed   
  using a subgroup  of $\Sp(2,q)=\SL(2,q)$ 
  sharply transitive on $\F_q^2-\{0\}\,$
  \cite{Pe,DH,CDS}.  It  is contained in the  
  non-symplectic spread of  $\F_q^4$
corresponding to the associated affine  irregular  nearfield plane.

  \vspace{-2pt}

\section{Concluding remarks}

The preceding examples makes it clear that there are rather few
known types of maximal symplectic partial spreads.  
There are amazingly few
known types  in odd characteristic, especially in view of the tables in \cite{CDFS,Grassl}.    
We mentioned a number of symplectic partial spreads whose maximality has yet to be decided.

We have not yet considered most inequivalence questions
for given dimension and field size.   
Suppose that $q$ is  even.
The  number  of   inequivalent orthogonal spreads in 
  $\O^+(4m ,q )$-spaces  is not bounded above by any polynomial in  $q^m$  
 \cite{KaW};
 these   produce inequivalent   maximal
     symplectic partial spreads
         \vspace{1pt}
 in  \propref{orthogonal spread}, \ref{nk}   
   and
\ref{project k}. 
In addition, 
there are at least 
   $q^{q^k}/q^{4k^2}$ 
   inequivalent examples
    in 
\thmref{Grassl example}(ii),  
$\binom{  q-1 }{s}/q^{30 }$ 
 for each pair $q,s$
    \vspace{2pt}
  in \thmref{orthovoids}, 
    $(  q-1 )^{q+1}/q^{30 }$ in
 \thmref{q2+q+1},
$\binom{  q+1 }{s}/q^{30 }$
\vspace{2pt}
 for each pair $q,s$
  in \thmref{smaller Grassl},
 $\binom{q+1 }{s}/q^{11 }$
 for each pair $q, s$
  in \thmref{conics} and 
$(q-2)(q-3)/6\log q$ in  \exref{3q-1}. 

\smallskip\smallskip
{\noindent\bf Acknowledgements.}
I am grateful to Markus Grassl for stimulating my interest in maximal symplectic partial spreads by pointing out the scarcity of   examples in 
dimension  $>4$.
This research was supported in part by a grant  from the Simons Foundation.


\appendix

\section{The desarguesian ovoid in $\O^+(8,q)$-space}
\label{appendix A}
In order to prove    \thmref{orthovoids}
we will consider a specific orthogonal ovoid
 in an  $\O^+(8,q)$-space with {$q>2$  even.} 
Let $F=\F_{q^3}\supset K=\F_q$, 
with trace map $T\col F\to K$ and norm $N\col F\to K$.
   Then $Q(a,\b,\g,d):=ad+T(\b\g)$ turns $V:=K\oplus F \oplus F \oplus K$ into 
an $\O^+(8,q)$-space. 

 The $q^3+1 $ points
$\<(0,0,0,1 ) \>$ and $\<(1,t,t^{q+q^2}, N(t )\>$,  $t\in F$,  
form  an ovoid
$\Omega$ on which  $G:=\SL(2,q^3)$ acting 3-transitively.  
In \cite[p.~1204]{Ka1}  this is called a {\em desarguesian ovoid}
(since it arises from  a desarguesian spread 
of an $\Sp(6,q) $-space  using
\lemref{orthogonal to symplectic}(ii) and  triality),
and it is observed that $G$   has exactly  two  orbits of singular points
of $V$,
one of which is~$\Omega$.
If $q>2$ and $p$ is any singular point not in $\Omega$, then 
$\<p^\perp\cap \Omega\>=p^\perp$ \cite[p.~1204]{Ka1}, as required in
\lemref{use orthogonal ovoid}. 
   
\notation
\label{dagger}\rm
Let $\pi\in F$ with $T(\pi)=0\ne T(\pi^{{1+q}} )$.   Use the 
nondegenerate symmetric $K$-bilinear  form $T(xy)$ on $F$ to see that
$\pi^{q} \notin\{t\in F\mid T(\pi t)=0\}= K+K\pi.$
   \vspace{-4pt}
 \Lemma
\label{non-perp}
If $p_1$ and $p_2$ are distinct singular points not in $\Omega ,$
then $|p_1^\perp \cap p_2^\perp \cap \Omega|\le  5q-5$.
\rm

\medskip
\proof
By  the transitivity of $G$ we may assume  
 that $p_1=\<(0,0,\pi,0)\>$ and 
$p_2=\<(a,\b,\g,d)\>$
for some  $a,\b,\g,d$.  We need to estimate the number of solutions $t$ to the equations 
$$T(t\pi)= 0=aN(t)+d+T(\b t^{q+q^2}   +\g t)  
$$
corresponding to points $\<(1,t,t^{q+q^2},t^{1+q+q^2})\>$.
By \eqref{dagger}
we can write $t=u+v\pi$ with $u,v\in K$.   Then the   second   equation
is    
$$
aN(u+v\pi)+d+T(\b[u+v\pi]^{q+q^2} +\g[u+v\pi])  = 0,
$$ 
which expands to 
\begin{equation}
\label{main equation2}
\begin{array}{lll}
a\{u^3+uv^2T(  \pi^{q+q^2}) + v^3N(\pi)\} +d 
  \vspace{2pt}
\\
\hspace{23pt}
 +~u^2T(\b ) +uvT(\b\pi) +v^2 T(\b \pi^{q+q^2}) +uT(\g) +v T(\g \pi)
= 0.
\end{array}
\end{equation}
For each $u$  this is a $K$-polynomial in $v$ of degree at most three,
and hence has  at most three roots $v\in K$ if it is not the zero polynomial.
Let $B$ be the number of   ``bad'' $u$ for which  this polynomial in $v$  is the zero
polynomial.  Then $|p_1^\perp \cap p_2^\perp \cap \Omega|\le  
(q-B)3+Bq +1$ (the last term  occurs since
$\<(0,0,0,1)\>$ may be in the intersection).  We will show that $B\le2$,
which produces  the bound in  the lemma.

The coefficients of our polynomial show that, for a ``bad'' $u$,
 we must have $a=0$,
then $T(\b \pi^{q+q^2}) =0$,$\,u T(\b\pi )+T(\g\pi)=0 $
and $ u^2T(\b ) + uT(\g)+ d =0$.  If $T(\b\pi )\ne0$ then there is 
one ``bad'' $u$, and if $ T(\b\pi )=T(\g\pi)=0 $ then  there are at most two
``bad'' $u$ unless $ T(\b ) = T(\g)= d =0$.

Thus, we must show that ${ T(\b \pi^{q+q^2}) =T(\b\pi )=T(\g\pi)= T(\b ) = T(\g)=0 }$
cannot all occur.~Since $  T(\b ) =T(\b\pi ) =0$, by
 \eqref{dagger}  we have $\b=x\pi$ with ${x\in K}$.~Then $0=T(\b \pi^{q+q^2}) =xT(N(\pi))$, so that $x=0$.
Similarly, $  T(\g ) =T(\g\pi ) =0$  implies that  $\g=y\pi$ with $y\in K$.  Now $p_2=\<(0,0,y\pi,0)\>=p_1$,
which is not the case.  \qed

\notation\rm
Let $\Omega_0\subset  \Omega$ consist of   $\<(0,0,0,1 ) \>$ and $\<(1,t,t^{q+q^2},t^{1+q+q^2})\>$,  $t\in K$.  There are $(q+1)^2$ singular points in $\Omega_0^\perp,$
all having the form $\<(0,\b,\g,0)\>$  with   $T(\b)=T(\g)=T(\b\g)=0 $.
The  sets $\Omega_0$  and $\Omega_0^\perp$  are   acted on by~a naturally embedded subgroup $G_0   =  \SL(2,q)$  of $G$ containing the transformations  %
\vspace{-2pt}
$$
\begin{array}{llll}
\hspace{-5.5pt}u_s\col (a,\b,\g,d)\mapsto 
(a,\b+sa,\g   +  as^2  +  \b^qs  +\b^{q^2}\!  s,d+as^3  +  T(\b)s^2  +  T(\g)s) 
\vspace{2pt}
\\
j\col  (a,\b,\g,d)\mapsto (d,\g, \b,a)   .
 \vspace{-2pt}
\end{array}%
$$
with $   s\in K$.  These act  on each of the $q+1$ lines 
$\<(0,\b,0, 0), (0,0,\b, 0) \>$ 
with $T(\b) =0\ne \b $
that partition the $(q+1)^2$ singular points
in $\Omega_0^\perp$,   sending 
\vspace{-1pt}
\begin{equation}
\label{generators}
\begin{array}{llll}
\hspace{-5.5pt}u_s\col (0,\b,\g,0)\mapsto 
(0,\b,\g   +  \b s  ,0)
\vspace{2pt}
\\
j\col (0,\b,\g,0)\mapsto (0,\g, \b,0)   . \hspace{185pt}
\end{array}%
\vspace{-1pt}
\end{equation}
\definition
An {\em ordinary} point is a 
singular point in  $\Omega_0^\perp$  of the form $\< (0,\b, \g ,0)\>$
such that either $\b=0$ and  $T(\g^{1+q})\ne 0$,
or    $T(\b^{1+q})\ne 0$  (recall
  that  $ T(\b)=T(\g)= T(\b\g) = 0 $).~Since any $\b\in F^*$
has characteristic polynomial $x^3+T(\b)x^2+T(\b^{1+q})+N(\b)$,
the ordinary requirement  can fail for some $\b,\g$ if 
and only if $q\equiv 1$ (mod 3). Moreover, 
if  $\b\in F-K$ then
${\b ^q \in \b K}\iff~{ \b^{q-1}\in K } \iff  \b^{(q-1,q^2+q+1)}\in K  \iff \b^{3}\in K 
\iff     T(\b^{1+q})= 0$.

For $\pi$ in
\eqref{dagger},  since $T((a\pi+ \pi^q)(a\pi + \pi^{q }) ^q )=(a^2+a+1)T(\pi^{1+q})$
the points of the line
 $\<(0,a\pi+ \pi^q ,0 ,0),   (0,0 , a\pi+  \pi^q ,0) \>,$ $a\in K,$  
  are ordinary if and only if 
 $a^2+a+1\ne 0$, so that 
 all   points   are ordinary if $q\equiv 2$ (mod 3),
 but   there are  two  lines of this form all of   whose points are not ordinary
 when    $q\equiv 1$ (mod 3).

The significance   of ordinary points is the following 
   \begin{lemma}
\label{move to ordinary}
If $p$ is an ordinary point then
\begin{itemize}
\item [\rm(i)]$p$ has the form  $ \<(0,0,\g, 0)\> $ with $T(\g)=0$ 
or $\<(0,\b,a \b, 0)\> $  with $T(\b)=0$ and $a\in K,$
and
\item  [\rm(ii)]$p^g=\<(0,0,\pi',0)\>$ for some $g\in G_0,$
where
$\pi'$ behaves as $\pi $ does  in \eqref{dagger}$:$ \
$T(\pi')=0\ne T(\pi'{}^{{1+q}})$.
\end{itemize}

   \end{lemma}
      
   \proof We may assume that  $p=\<(0,\b,\g, 0)\> $  with $\b\ne0$.
   
   (i)
   Since $p$ is ordinary,  $\b^q\notin K\b $,
   so that    $\b$ and $\b^q$ span $\ker T $.
   Write $\g=k\b+b\b^q$ with $k,b\in K$.  Then 
   $0=T(\b\g)=bT(\b^{{1+q}})$ implies that $b=0$.
   
   (ii)
By \eqref{generators},  $p^{u_k j}=\<(0,0,\b,0)\>$ behaves as stated.
\qed

   \begin{lemma}
\label{new intersection sizes}
If $p_1,p_2,$   and  $p_3$
 are pairwise non-perpendicular ordinary points$,$
then
\begin{itemize} 
\item[\rm(i)] $ |p_1^\perp \cap p_2^\perp \cap \Omega | =2q$  and
\item[\rm(ii)]  
 $ |p_1^\perp \cap p_2^\perp \cap p_3^\perp  \cap \Omega | =q+2. $
\end{itemize}
\end{lemma}
   
   \proof
By \lemref{move to ordinary}(ii)  we may assume that $p_1$ has the form 
$\<(0,0,\pi,0)\>$ and 
$p_2 = \<(0,\b,\g,0)\>$, where $T(\b)=T(\g)=T(\b\g)=0$.
Also $ T(\b\pi)\ne0$
since $p_1$ and $p_2$ are not perpendicular.
All $(0,0,0,1) $  and  $(1,t,t^{q+q^2},N(t)),t \in K$,
 are in each of the stated intersections, 
so we will focus on vectors $(1,t,t^{q+q^2},N(t))$
with $t=u+v\pi \notin K$   in the   intersections.

(i) Here \eqref{main equation2} states that 
\begin{equation}
\label{to massage}
uv T(\b\pi ) +v^2T(\b\pi^{q+q^2})+vT(\g\pi)=0.
\end{equation}
Since  $T(\b\pi)\ne0,$  each $v\ne0$ 
 determines a unique $u$.  This argument reverses: the 
 intersection size  is  $ (q+1)+(q-1)$.

 Before continuing we  massage  \eqref{to massage}. 
By   \lemref {move to ordinary}(i), 
 $\g=k\b$  for some $k\in K$. 
  Since $\dim \ker T =2$  we can write $\b=x\pi +y\pi^q$ with $x,y\in K$.
 Since $0\ne T(\b\pi) =yT(\pi^{1+q} )$ we have
 $y\ne0$ and  $\b\in ((x/y)\pi + \pi^q)K$.  We may assume that 
$\b = a\pi + \pi^q $ with $a\in K$.   Then
\begin{equation}
\label{p2}
p_2 = \<(0,a\pi + \pi^q ,k (a\pi + \pi^q ),0)\>,
\end{equation}
so that $T(\b\pi)=T(\pi^{{1+q}})$  and   \eqref{to massage}   becomes
\begin{equation}
\label{mn eq}
uT(\pi^{{1+q}}) + v[aN(\pi)+ T(\pi^{2q+q^2})]+kT(\pi^{{1+q}}) =0.
\end{equation}

(ii)    We may assume that $p_3=\<(0,\b',\g',0\>)$ with 
$\g'=k'\b'$ and $\b' =a'\pi + \pi^q $  for some $k',a'\in K$.
Then  $(a+a')(k+k')T(\pi\pi^q)=  T(\b\g'+\g\b')\ne0$.
The two versions  of  \eqref{mn eq} imply that
\begin{equation}
\label{mn}
v=\frac{k+k'}{a+a'}\frac{T(\pi^{{1+q}})}{N(\pi)}, \ \ 
u=k+  \frac{k+k'}{a+a'} \Big(a+\frac{T(\pi^{2q+q^2}) }{ N(\pi)}\Big),
\end{equation}
which proves (ii).
 \qed

\begin{Example}
\label{examples}
{\rm(i)}    If $ {\mathcal P}\subseteq
   \{\<(0,0,\pi,0)\>,   \<(0,a\pi+ \pi^q ,a^2\pi+a\pi^q ,0)\> \mid {a \in K,}~
   a^2+a+1\ne0\},$ then 
   \vspace{-6pt}
 $$\displaystyle  \Big |\bigcap_{p\in {\mathcal P} } p ^\perp   \cap \Omega \Big | =
   \begin{cases} 
   {q^2+1} &\mbox{if}   \  \  |{\mathcal P} | =1 \\
   {2q}	&\mbox{if}      \  \  |{\mathcal P} |=2 \\
   {q+2}	&\mbox{if}   \ \   |{\mathcal P} |\ge 3  .
 \end{cases}
  $$
{\rm(ii)}
   If ${\mathcal P} \subseteq  
   \{\<(0,0,\pi,0)\>, \ 
   \<(0,  \pi^q ,0 ,0)\>,\ 
   \<(0, \pi+ \pi^q , \pi+ \pi^q ,0)\>, 
   \<(0,a\pi+ \pi^q ,a^3\pi+a^2\pi^q ,0)\>
  \} $  for  an arbitrary   $a\in K-\{0,1\} $ 
  such that $ a^2+a+1\ne0,$ then    
   \vspace{-2pt}
 $$\displaystyle \Big |\bigcap_{p\in {\mathcal P} } p ^\perp   \cap \Omega \Big | =
   \begin{cases} 
   {q^2+1} &\mbox{if}      \  \  |{\mathcal P} |=1 \\
   {2q}	&\mbox{if}      \  \  |{\mathcal P} |=2 \\
   {q+2}	&\mbox{if}   \ \   |{\mathcal P} |=3    \\
   {q+1}	&\mbox{if}   \ \   |{\mathcal P} |=4 .
 \end{cases}   \vspace{-2pt}
  $$
\end{Example}
\proof
All of the stated points are ordinary.

(i)  In
  \eqref{p2},  $k  =  a$ for all listed points other than 
$\<(0,0,\pi,0)\>$.~By   \eqref{mn},
$t =  
\frac{T(\pi^{2q+q^2}) }{ N(\pi)\raisebox{1.5ex}{\hspace{-1pt}}\raisebox{-.4ex}{\hspace{-1pt}}}  
+  \frac{T(\pi^{{1+q}})}{ N(\pi)\raisebox{1.5ex}{\hspace{-1pt}}\raisebox{-.4ex}{\hspace{-1pt}}}\pi  
$
is in every intersection (which is easily 
checked directly); so is $\Omega_0$, so that
  every intersection~has size   $\geq q+2$. Since any intersection of~three sets $p^\perp\cap \Omega$ has  size $q+2$ (by  \lemref{new intersection sizes}(ii)),  so does 
  any intersection of at least four such~sets.%
  \medskip

(ii)  The last three of these four ordinary points correspond to the pairs 
$(a,k)=(0,0),\,(1,1), \,(a,a^2)$ in \eqref{p2}.  
Then   \eqref{mn} 
and  different  3-sets in 
$\mathcal P$ produce     different values of $v$, so that 
$ { |\bigcap_{p\in {\mathcal P}} p ^\perp   \cap \Omega  | }=q+1$ if 
$| {\mathcal P}|=4$.~The remaining sizes  are given in \lemref{new intersection sizes}.\qed
 
\appendix
\renewcommand{\thesection}{B}

\section{Suzuki-Tits  ovoids: background}
\label{Background on Suzuki-Tits  ovoids}
 
We will need   information  
concerning     a    Suzuki-Tits  ovoid $\Omega$ in an
 $\O(5,q)$-space  $U$  with radical $r$, where $q=2^{2e+1}$.
 The standard view of these ovoids is in symplectic space.
 For our purposes, the  view from an $\O(5,q)$-space     has advantages, such as 
 lying in an $\O^+(8,q)$-space. 
 
 Let $\bar\Omega$ denote a standard Suzuki-Tits ovoid in the symplectic 4-space $U/r$\ \cite{Ti2}.
If  $\<x,r\>/r\in \bar\Omega$ then the line $\<x,r\>$
has a unique singular point.   Thus, there is a set 
$\Omega $ of $q^2+1$ singular points of $U$ that projects onto 
$\bar\Omega$.    
The group   $ \Sz(q) $ lifts from a subgroup of $\Sp(4,q)$
to a  group  $G < \O(5,q)$ preserving  $ \Omega $.    
See  \cite[Sec.~6.4]{Demb}  for information concerning the inversive plane 
$\I(\Omega)$ produced by  $ \Omega$.

We will assume that  $q>2$.  Then   $U=\<\Omega\>$
since $G$ does not act on an $\O^{\pm}(4,q)$-space.
(If $q=2$ then $\Omega$ spans an 
$\O^-(4,2)$-space.)
 
 \Lemma
 \label{hyperplane intersections}
 Every hyperplane meets $\Omega$.  
 More precisely$,$     there are exactly five types of hyperplanes $H$ 
 of $U$$:$
  \begin{itemize}
 \item[\rm (i)] Tangent hyperplanes $p^\perp$ for $p\in \Omega,$
 with $r\in H$ and $H\cap \Omega=\{p\};$
  \item[\rm (ii)]  Secant hyperplanes $x^\perp=H$  containing $r,$ 
  where   $x$ is a singular point$,$  $x^\perp\cap \Omega$ is a circle of  $\I(\Omega)$   and $\<x^\perp\cap \Omega\>=x^\perp;$
  \item[\rm (iii)] $\O^-(4,q)$-hyperplanes  for which
 $H\cap \Omega$ is an orbit of a cyclic subgroup of $G$ of order
 $|H\cap \Omega|=q\pm\sqrt{2q}+1$ acting   irreducibly on $U/ r;$
  and
   \item[\rm (iv)]  $\O^+(4,q)$-hyperplanes for which  
 $ H\cap \Omega  $
  contains an orbit of a cyclic subgroup of $G$ of order
 $|H\cap \Omega|-2=q -1$ that fixes two  points of  $H\cap \Omega .$  Moreover$, $
 $ H\cap \Omega  $ is not one of the circles in {\rm  (ii)}.
  \end{itemize}
\rm  
\smallskip

\proof
(i) Projecting mod $r$ shows that each point of $\Omega$ behaves as
stated.

\smallskip
(ii) If $x$ is a singular point not in $\Omega$
then each of the $q+1$ t.s. lines on $x$ meets $\Omega$ 
since $\Omega$  is an ovoid, so that  $|x^\perp\cap \Omega|=q+1$.
Also,  $\dim\<x^\perp\cap \Omega\>=4$,  as 
 otherwise its set of  singular points would    project into  a plane of $U/r$,
  and hence be contained in  a conic,
which is not the case  since $q>2$   
\cite[pp.~51-52]{Ti3}. 
Since $\<x^\perp\cap \Omega\>$ lies in the 4-space $x^\perp$,
these subspaces   coincide

\smallskip
 (iii)   This is  
  \cite[Theorem~1(a)]{Bagchi-Sastry2}. 
  
  \smallskip
  (iv)    The set of singular points of $H$ is partitioned by $q+1$ t.s.~lines, and each
  t.s.~line of $U$ meets $\Omega$ since $\Omega$  is an ovoid.
  Thus, $|H\cap \Omega|=q+1$.

We use  the orbits of $G$ to find $G_H$.  
There are  exactly two point-orbits on $U/r$:
$\bar\Omega$  and the remaining $q(q^2+1)$ points.
There is a subgroup of $G$ of order $q-1$ that fixes four points of $U/r$
and induces all scalars on each of these 1-spaces  \cite[p.~183]{HB}.
Since each  line containing  $r$ 
has a unique singular point,  the two point-orbits on 
 $U/r$ produce four  point-orbits on  $U-\{r\}$.

Since  $G$ has five point-orbits it  also has  five hyperplane-orbits, so that 
    all $q^2(q^2+1)/2$ hyperplanes $H$ in (iv) lie in an orbit. 
Then $|G_H|= |G| /[q^2(q^2+1)/2]=2(q-1)$,
so that  
 $G_H$ is dihedral of order $2(q-1)$, 
with orbits of size $2$ and $q- 1$ on  $\Omega$ \cite[Theorem~9]{Suz}.

For the final assertion, if $H\cap \Omega$ lies in two hyperplanes then it is in a plane, and hence is a conic, which is not the case
\cite[pp.~51-52]{Ti3}. 
\qed  

\end{document}